\documentclass[11pt]{amsart}
\usepackage{palatino}
\usepackage{amsfonts}
\usepackage{amssymb}
\usepackage{amscd}

\theoremstyle{plain}

\newtheorem{theorem}{Theorem}[section]
\newtheorem{proposition}[theorem]{Proposition}
\newtheorem{corollary}[theorem]{Corollary}
\newtheorem{lemma}[theorem]{Lemma}
\theoremstyle{definition}
\newtheorem{definition}[theorem]{Definition}
\newtheorem{remark}[theorem]{Remark}

\newtheorem{conjecture/question}[theorem]{Conjecture/Question}

\newtheorem{remark/definition}[theorem]{Remark/Definition}
\newtheorem{terminology/notation}[theorem]{Terminology/Notation}

\newcommand{\marginlabel}[1]%
  {\mbox{}\marginpar{\raggedleft\hspace{0pt}\bfseries\sf#1}}

\def\PP{{\textbf P}}
\def\OO{\mathcal{O}}
\def\cN{\mathcal{N}}

\def\cB{\mathcal{B}}
\def\cA{\mathcal{A}}

\def\E{\mathcal{E}}
\def\G{\mathcal{G}}

\def\L{\mathcal{L}}

\def\cM{\mathcal{M}}

\def\cU{\mathcal{U}}

\def\H{\mathcal{H}}
\def\Pic0{{\rm Pic}^0(X)}

\def\mm{\overline{\mathcal{M}}}
\def\kk{\overline{\mathcal{K}}}
\def\zz{\overline{\mathcal{Z}}}
\def\KK{\overline{\mathcal{K}}}
\def\GG{{\textbf G}}
\def\GisP{{\overline{\mathcal{GP}}_{g,d}^r}}

\theoremstyle{remark}

\begin{document}

\title{\bf Koszul divisors on moduli spaces of curves}

\author[G. Farkas]{Gavril Farkas}
\address{Department of Mathematics, University of Texas,
Austin, TX 78712} \email{{\tt gfarkas@math.utexas.edu}}
\thanks{Research  partially supported by an Alfred P. Sloan Fellowship,  the NSF Grants DMS-0450670 and DMS-0500747
and  a 2006 Texas Summer Research Assignment}

\maketitle

\section{Introduction}

In this paper we describe a general method of constructing special
effective divisors on various moduli spaces using the syzygies of
the parametrized objects. The method can be applied to
a wide range of moduli problems with the property that the coarse
moduli space has canonical singularities hence pluricanonical forms
 extend over any desingularization of the moduli space. Here we
treat the case of the moduli stacks $\mm_{g, n}$ and we develop the
intersection theory machinery necessary to understand the
compactification and compute the class of these Koszul divisors. Our main result
(Theorem 1.1) provides the first infinite sequence of actual (as opposed to virtual) counterexamples to
the Harris-Morrison Slope Conjecture and encodes in a single formula virtually all known divisor
class calculations on $\mm_g$.

The idea of using geometric divisors to study the  geometry of
 a moduli space can be traced back to Harris and Mumford  (cf.  \cite{HM})
 who, in the course of their proof  that $\mm_g$ is of general type for
 odd genus $g=2k-1\geq 25$, studied the \emph{Hurwitz} divisor
 $\cM_{g, k}^1:=\{[C]\in \cM_{g}: \exists \mbox{ } C\stackrel{k:1}\rightarrow  \PP^1\}$
  consisting of curves with a pencil $\mathfrak g^1_k$. By computing the class of $\mm_{g, k}^1$ and comparing it
  in to  $K_{\mm_g}$, they showed that when $g\geq 25$,  the canonical class
  is a combination with positive coefficients of $[\mm_{g, k}^1]$, the Hodge class $\lambda$ (which is big and nef)
 and various boundary classes. Later, numerous other divisor class calculations were carried out.
 Eisenbud and Harris considered the \emph{Petri divisors} on $\mm_{2k-2}$ consisting of curves $C$ of genus
 $2k-2$ having a pencil $A\in W^1_k(C)$ which violates the Petri
 Theorem,
  which then they used to show that $\mm_g$ is of general type for even $g\geq 24$ (cf.
  \cite{EH3}).
Logan introduced  \emph{pointed} Brill-Noether  divisors on $\mm_{g,
g}$ consisting of curves $[C, x_1, \ldots , x_g]\in \mm_{g, g}$ with
the property that $h^0(C, \OO_C(x_1+\cdots +x_g))\geq 2$
 (cf. \cite{Log}) and used them to determine the Kodaira type of $\mm_{g, n}$ for various $g$ and $n$.

More recently, in \cite{FP} in our work on the Harris-Morrison Slope
Conjecture, we  reinterpreted the condition that a curve $[C]\in
\cM_{10}$ lie on a $K3$ surface as saying that there exists a linear
system $L=K_C(-\mathfrak g^1_6) \in W^4_{12}(C)$ such that the
embedded curve $C\stackrel{|L|}\hookrightarrow \PP^4$ is not
projectively normal. Using this description we computed the class of
the compactification of the divisor $\mathcal{K}_{10}$ of curves
with this property and showed that $s(\KK_{10})=7$, thus
contradicting the Slope Conjecture. In \cite{F3} we generalized this
construction to cover all cases $g=6i+10$ and we obtained a
(sometimes virtual) Hurwitz type divisor on $\mm_{6i+10}$ defined in
terms of linear series $\mathfrak g^{3i+4}_{9i+12}= K_C(-\mathfrak
g^1_{3p+6})$ residual to a pencil of minimal degree. This locus,
when a divisor, always has slope $<6+12/(g+1)$ thus violating the
Slope Conjecture (see \cite{HMo} and \cite{FP} for background on the
effective cone of $\mm_g$ and for the significance of the
Harris-Morrison Conjecture). Around the same, Khosla provided a
different type of example of a divisor on $\mm_g$ having
exceptionally small slope (cf. \cite{Kh}): on $\mm_{21}$ the closure
of the locus of curves $[C]\in \mm_{21}$ possessing an embedding
$C\hookrightarrow \PP^6$ given by a $\mathfrak g^6_{24}$ such that
$C$ lies on a quadric, is a divisor whose slope is less than the
slope of the Harris-Mumford divisor $\mm_{21, 11}^1$.

The aims of this paper are (1) to give a unified framework for doing
divisor class calculation on $\mm_{g, n}$ and (2) to provide (empirical)
evidence that syzygy divisors may be the answer to the riddle: Given
a moduli space, what is the most intrinsic, most natural and from
the point of view of birational geometry, most useful effective
divisor on it? We prove that virtually all interesting known
divisors on $\mm_g$ (the Harris-Mumford divisor, the Petri divisor
and all known counterexamples to the Harris-Morrison Conjecture) can
be treated in a unified way and are particular instances of a single
syzygy type construction. In \cite{F4} we shall further illustrate
this ideology by studying moduli spaces of curves with various level
structures from the point of view of syzygies.

We fix integers $i\geq 0$ and $s\geq 1$ and set $r:=2s+si+i,
g:=rs+s$ and $d:=rs+r$. We denote by $\mathfrak G^r_d$ the stack
parametrizing pairs $[C, L]$ with $[C]\in \cM_g$ and $L\in W^r_d(C)$
and denote by $\sigma:\mathfrak G^r_d\rightarrow \cM_g$ the natural
projection. Since $\rho(g, r, d)=0$,  by general Brill-Noether
theory, the general curve of genus $g$ will have finitely many
$\mathfrak g^r_d$'s and there exists a unique irreducible component
of $\mathfrak G^r_d$ which maps onto $\cM_g$.

We denote by $K_{i, j}(C, L)$ the $(i, j)$-th Koszul cohomology
group of the pair $[C, L]\in \mathfrak G^r_d$ and define a
stratification of $\mathfrak G^r_d$ with strata $\mathcal{U}_{g,
i}:=\{(C, L)\in \mathfrak G^r_d: K_{i, 2}(C, L)\neq 0\}$. We then
set $\mathcal{Z}_{g, i}:=\sigma_*(\mathcal{U}_{g, i})$.

\begin{theorem}\label{div} If $\sigma:\widetilde{\mathfrak G}^r_d\rightarrow
\widetilde{\cM}_g$ is the compactification of $\mathfrak G^r_d$
given by limit linear series on tree-like curves, then there exists
a natural morphism between torsion free sheaves of the same rank
$\phi: \cA\rightarrow \cB$ over $\widetilde{\mathfrak G}^r_d$ such
that $\overline{\mathcal{Z}}_{g, i}$ is the image of the degeneracy
locus of $\phi$. The class of the pushforward to $\widetilde{\cM}_g$
of the virtual degeneracy locus of $\phi$ is given by
$$\sigma_*(c_1(\cB-\cA))\equiv a\lambda-b_0\ \delta_0-b_1\ \delta_1-\cdots -b_{[g/2]}\ \delta_{[g/2]},
  $$
where $a, b_0, \ldots, b_{[g/2]}$ are explicitly given coefficients
such that $b_1=12b_0-a$ and
$$s\bigl(\sigma_*(c_1(\cB-\cA))\bigr)=\frac{a}{b_0}=6\frac{f(s,
i)}{(i+2)\ s \ h(s, i)}, \mbox{ with} $$ \begin{center} $f(s,
i)=(i^4+8i^3+24i^2+32i+16)s^7+(i^4+4i^3-16i-16)s^6-(i^4+7i^3+13i^2-12)s^5-
(i^4+2i^3+i^2+14i+24)s^4
+(2i^3+2i^2-6i-4)s^3+(i^3+17i^2+50i+41)s^2+(7i^2+18i+9)s+2i+2$
\end{center}
and \noindent
\begin{center}
$h(s,i)=(i^3+6i^2+12i+8)s^6+(i^3+2i^2-4i-8)s^5-(i^3+7i^2+11i+2)s^4-\newline
-(i^3-5i)s^3+(4i^2+5i+1)s^2+ (i^2+7i+11)s+4i+2.$
\end{center}
Furthermore, we have that $6<\frac{a}{b_0}<6+\frac{12}{g+1}$
whenever $s\geq 2$. If the morphism $\phi$ is generically
non-degenerate, then $\overline{\mathcal{Z}}_{g, i}$ is a divisor on
$\mm_g$ which gives a counterexample to the Slope Conjecture for
$g=s(2s+si+i+1)$.
\end{theorem}

For a precise definition of the partial compactification
$\widetilde{\cM}_g\subset \mm_g$ of $\cM_g$ we refer to Section 2.
Since $\mathrm{codim}(\mm_g-\widetilde{\cM}_g, \mm_g)\geq 2$, it
makes no difference whether the computation of
$[\overline{\mathcal{Z}}_{g, i}]$ is carried out over
$\widetilde{\cM}_g$ or $\mm_g$. Despite its complicated appearance,
the slope computed in Theorem \ref{div} encodes a surprising amount
of information about $\mm_g$. In particular, for suitable choices of
$s$ and $i$ it specializes to the divisor class calculations carried
out in \cite{HM}, \cite{EH3}, \cite{Kh}, \cite{FP} and \cite{F3}
which were originally obtained using a variety of ad hoc techniques.
The first interesting case is $s=1, g=2i+3$ when $\mathfrak
g^r_d=\mathfrak g^{g-1}_{2g-2}=K_C$ (the canonical bundle is the
only $\mathfrak g^{g-1}_{2g-2}$ on a curve of genus $g$). We can
relate the locus $\mathcal{Z}_{2i+3, i}$ to more classical loci in
$\cM_{2i+3}$ using Green's Conjecture which predicts that for any
smooth curve $C$ one has the equivalence $K_{l, 2}(C,
K_C)=0\Leftrightarrow l<\mbox{Cliff}(C)$. Although Green's
Conjecture for arbitrary curves is still open, Voisin proved it for
generic curves of given gonality (cf. \cite{V1}, \cite{V2}). In our
case this gives a set-theoretic identification between
$\mathcal{Z}_{2i+3, i}$ and the locus $\mathcal{M}_{2i+3, i+2}^1$ of
$(i+2)$-gonal curves. Thus Theorem \ref{div} provides a new way of
calculating the class of the compactification of the Brill-Noether
divisor first computed by Harris and Mumford (cf. \cite{HM}):

\begin{corollary}\label{bn}

The slope of the Harris-Mumford divisor $\mm_{2i+3, i+2}^1$ on $\mm_{2i+3}$  consisting
of curves which  cover  $\PP^1$ with degree $\leq i+2$ is given by the
formula $$s(\mm_{2i+3, i+2}^1)=\frac{6(i+3)}{i+2}=6+\frac{12}{g+1}.$$
\end{corollary}

For $s=2$ and $g=6i+10$ (that is, in the case $h^1(L)=2$ when
$\mathfrak G^r_d$ is isomorphic to a Hurwitz stack parameterizing
covers of $\PP^1$), we recover the main  result from \cite{F3}:
\begin{corollary}
The slope of the divisor $\overline{\mathcal{Z}}_{6i+10, i}$ on $\mm_{6i+10}$ consisting
of curves possessing a pencil $\mathfrak g^1_{3i+6}$ such that if $L=K_C(-\mathfrak g^1_{3i+6})\in W^{3i+4}_{9i+12}(C)$
denotes the residual linear system, then  $C \stackrel{|L|}\hookrightarrow \PP^{3i+4}$ fails to satisfy the Green-Lazarsfeld property $(N_i)$,
is given by the formula:
$$s(\overline{\mathcal{Z}}_{6i+10,
i})=\frac{3(4i+7)(6i^2+19i+12)}{(12i^2+31i+18)(i+2)}.$$
\end{corollary}

In the case $i=0$ we have complete results in the sense that (1) we
show that $\zz_{g, 0}$ is an actual divisor on $\mm_g$ and (2) we
can compute the entire class  $[\zz_{g, 0}]$ rather than the
$\lambda, \delta_0$ and $\delta_1$ coefficients. In particular we
show that $b_j\geq b_0$ for $j\geq 1$, hence the slope of $\zz_{g, 0}$
is always
computed by the $\lambda$ and $\delta_0$ coefficients.
\begin{theorem}\label{khosla}
For $g=s(2s+1), r=2s, d=2s(s+1)$ the slope of the virtual class of
the locus of those curves $[C]\in \mm_g$ for which there exists $L\in
W^r_d(C)$ such that the embedded curve $C\stackrel{|L|}\hookrightarrow \PP^r$ sits on a
quadric hypersurface, is
$$s(\overline{\mathcal{Z}}_{s(2s+1),
0})=\frac{a}{b_0}=\frac{3(16s^7-16s^6+12s^5-24s^4-4s^3+41s^2+9s+2)}{s
(8s^6-8s^5-2s^4+s^2+11s+2)}.$$
\end{theorem}

Note that this locus has been first considered by D. Khosla who,
using a different approach, was able to compute the coefficients $a$
and $b_0$ (cf. \cite{Kh}).  Showing that the degeneration loci
$\mathcal{Z}_{g, i}$ are actual divisors on $\cM_g$ can be very
difficult in practice (for instance, the statement that
$\mathcal{Z}_{2i+3, i}$ is a divisor on $\cM_{2i+3}$ is essentially
Green's Conjecture for a generic curve of odd genus). Apart from the
case $s=1$ (settled by Voisin in \cite{V2}), the only cases where it
was previously known that $\mathcal{Z}_{g, i}$ is an actual divisor
were $s=2, i=0$ (cf. \cite{FP}, this being the $K3$ divisor on
$\cM_{10}$), $s=2, i=1, 2$ (cf. \cite{F3}) and when $s=3, i=0$ (cf.
\cite{Kh}) - these last three cases having been settled using
Macaulay. Here we show that for $i=0$ the degeneracy loci
$\mathcal{Z}_{g, 0}$ are honest divisors on $\mm_g$, that is, the
map $\phi:\cA\rightarrow \cB$ described in Theorem \ref{div} is
generically non-degenerate. This provides the first infinite
sequence of actual (as opposed to virtual) counterexamples to the
Harris-Morrison Slope Conjecture \cite{HMo}:

\begin{theorem}\label{maxrank}
For an integer $s\geq 2$ we set $r:=2s, d:=2s(s+1)$ and
$g:=s(2s+1)$. Then  $\phi:\cA\rightarrow \cB$ is a generically
non-degenerate map between vector bundles over $\widetilde{\mathfrak
G}^r_d$ having the same rank and its degeneracy locus
$$\mathcal{Z}_{g, 0}:=\{[C]\in \cM_g: \exists L\in W^r_d(C) \mbox{ such that } C\stackrel{|L|}\hookrightarrow \PP^r
\mbox{ is not projectively normal}\}$$ is a divisor on $\mm_g$ of
slope
$$s(\overline{\mathcal{Z}}_{g,
0})=\frac{3(16s^7-16s^6+12s^5-24s^4-4s^3+41s^2+9s+2)}{s(8s^6-8s^5-2s^4+s^2+11s+2)}$$
contradicting the Slope Conjecture.
\end{theorem}

As an application of the techniques developed for proving Theorem
\ref{div} we compute the class of the Gieseker-Petri divisors on
$\mm_g$. Recall that Petri's Theorem asserts that for a
\emph{general} curve $[C]\in \cM_g$ and for an \emph{arbitrary} line
bundle $L$ on $C$, the multiplication map $$\mu_0(L): H^0(L)\otimes
H^0(K_C\otimes L^{\vee})\rightarrow H^0(K_C)$$ is injective (see
\cite{EH4} and \cite{Laz} for two very different, relatively short
proofs). The map $\mu_0(L)$ governs the deformation theory of
sections of the line bundle $L$. It is well-known that $G^r_d(C)$ is
smooth of expected dimension $\rho(g, r, d)$ at a point $[L]\in
W^r_d(C)$ if and only if $\mu_0(L)$ is injective. The locus in
$\cM_g$ where the Petri Theorem fails breaks up into numerous
components and its geometry is still quite mysterious (see
\cite{F2}, \cite{EH3}). For integers $r, s\geq 1$ we set again
$d:=rs+r$ and $g:=rs+s$, so that $\rho(g, r, d)=0$. Like in
\cite{F2} we define the Gieseker-Petri locus
$$\mathcal{GP}^r_{g, d}:=\{[C]\in \cM_g: \exists L\in W^r_d(C) \mbox{ such that }
\mu_0(L) \mbox{ is not injective }\}.$$

\begin{theorem}\label{gp}
For $d=rs+r$ and $g=rs+s$, the class of the Gieseker-Petri divisor
in $\mm_g$ is given by the formula:
$$\overline{\mathcal{GP}}^r_{g, d}\equiv \frac{c_r (s-1)r}{(r+s+1)(rs+s-2)(rs+s-1)}(a \lambda-b_0 \delta_0-b_1
\delta_1 -\sum_{j=2}^{[g/2]} b_j \delta_j), $$ where $c_r$ is an
explicitly given constant defined in Lemma \ref{vandermonde},
$$a=r^2 s^2 (4s+r+rs+10)+s^2(5rs+24r+2s+15)+21s+26rs+7r^2s+2r+2, $$
$$b_0=\frac{s(s+1)(r+1)(r+2)(rs+s+4)}{6}, $$
$$b_1=(rs+s-1)(3rs^2+2s^2+r^2s^2+7s+6rs+r^2s+2r+2),$$
and $ b_j\geq b_1$ for $ j\geq 2.$ In particular we have the
following expression for the slope:
$$s(\overline{\mathcal{GP}}_{g, d}^r)=6+ \frac{12}{g+1}+\frac{6(s+r+1)(rs+s-2)(rs+s-1)}{s(s+1)(r+1)(r+2)(rs+s+4)(rs+s+1)}.$$
\end{theorem}

Theorem \ref{gp} shows that the Gieseker-Petri divisors satisfy the
Slope Conjecture, that is, $s(\GisP)\geq 6+12/(g+1)$. This is
consistent with Proposition 2.2 from \cite{FP} stating that any
effective divisor on $\mm_g$ violating the Slope Conjecture would
have to contain the locus $\kk_g\subset \mm_g$ of curves lying on
$K3$ surfaces and with Lazarsfeld's result (cf. \cite{Laz}) that a
general $[C]\in \kk_g$ satisfies Petri's Theorem. For $s=2$, Theorem
\ref{gp} specializes to Eisenbud and Harris's computation originally
used to show that $\mm_g$ is of general type for large even genus
(cf. \cite{EH3}, Theorem 2):
\begin{corollary} \label{ehpetri}
For $g=2r+2$, the Gieseker-Petri divisor $\mathcal{GP}^r_{2r+2, 3r}$
can be interpreted as the branch locus of the generically finite map
$\sigma: \mathfrak{G}^r_{3r}\rightarrow \cM_{2r+2}$ from the Hurwitz
stack $\mathfrak{G}^r_{3r}=\mathfrak{G}^1_{r+2}$ of covers of degree
$r+2$ and one has the following expression for its class:
$$\overline{\mathcal{GP}}_{2r+2, 3r}^r\equiv
c_r\bigl(\frac{6r^2+25r+20}{2r+1}\lambda-\frac{(r+1)(r+2)}{2r+1}\delta_0-(3r+4)
\delta_1-\sum_{j=2}^{r+1} b_j \delta_j\bigr),$$ where $b_j>1$ for
$j\geq 2$.
\end{corollary}

In Section 4 we describe five different ways of constructing Koszul
divisors on $\mm_{g, n}$. The direct analogue of Theorem \ref{div}
in the pointed case is the following statement:
\begin{theorem}\label{pointedsyzygies}
Fix positive integers $g$ and $i$ such that
$$n:=\frac{2g+i+1}{2}+\frac{\sqrt{(i+1)^2+4ig+8g}}{2}$$ is an integer.
Then the locus $$\mathfrak{Syz}_{g, n}:=\{[C, x_1, \ldots, x_n]\in
\cM_{g, n}: K_{i, 2}\bigl(C, \OO_C(x_1+\cdots +x_n)\bigr)\neq 0\}$$
is a divisor on $\cM_{g, n}$, and the class of its compactification
is given by the formula:
$$\overline{\mathfrak{Syz}}_{g, n}\equiv \frac{1}{n-g-i}{n-g-1 \choose i}\Bigl(-(n+g-1)\lambda +(3g-n+i+1)\sum_{j=1}^n \psi_j-0\cdot \delta_{irr}-\sum_{j, t\geq 0}\sum_{|S|=t} b_{j:t} \delta_{j:S}\Bigr),$$
where $b_{j:t}>1$ are explicitly determined coefficients.
\end{theorem}

Another infinite sequence of interesting divisors on $\mm_{g, n}$
can be obtained by using the \emph{Gaussian-Wahl} map associated to
a line bundle on a curve. Recall that if $L$ is a line bundle on a
curve $C$, the Wahl map $$\psi_L:\wedge^2 H^0(C, L)\rightarrow
H^0(C, K_C\otimes L^{\otimes 2})$$ is defined by the formula
$\psi_L(f\wedge g):=f\cdot dg-g\cdot df$. The Gaussian $\psi_L$ measures deformations of the cone over the curve $C$ embedded in projective space by the linear system $|L|$ and it is  known that if $C$ lies on a $K3$ surface then the Wahl map $\psi_{K_C}$ cannot be surjective (cf. \cite{Wa}). Furthermore, the
divisor $\zz_{10, 0}$ on $\mm_{10}$ can be viewed as the global
degeneracy locus corresponding to the Wahl map for canonical
curve of genus $10$ (see \cite{FP} for details and further references). If $[C,
x_1, \ldots, x_n]\in \cM_{g, n}$, we set $\Gamma:=x_1+\cdots+x_n\in
C_n$ for the divisor of marked points.
\begin{theorem}\label{gauss}
Fix an integer $g$ such that  $$n:=\frac{2g+3+\sqrt{24g+1}}{2}$$ is
an integer. Then the locus $$\mathfrak{Wahl}_{g, n}:=\{[C, x_1,
\ldots, x_n]\in \cM_{g,n}: \psi_{\Gamma}: \wedge^2
H^0(\OO_C(\Gamma))\rightarrow H^0(K_C\otimes \OO_C(2\Gamma)) \mbox{
is degenerate}\}$$ is a divisor on $\cM_{g, n}$ and its
compactification has the following class:
$$\overline{\mathfrak{Wahl}}_{g,n}=-(n-g-1)\lambda+(n-g-1)\sum_{j=1}^n \psi_j-\delta_{irr}-\sum_{j,t\geq 0}b_{j:t}\sum_{|S|=t} \delta_{j:S},$$
where $b_{j:t}>1$ are explicitly determined coefficients.
\end{theorem}
Note that although the divisors $\overline{\mathfrak{Syz}}_{g, n}$
and $\overline{\mathfrak{Wahl}}_{g, n}$ live on $\mm_{g, n}$'s for
some very particular choices of $n$, using the forgetful and
clutching maps $\mm_{g, n}\rightarrow \mm_{g, n-1}$ and $\mm_{i,
n_1}\times \mm_{g-i, n_2}\rightarrow \mm_{g, n_1+n_2}$ one
immediately has explicit Koszul divisors on $\mm_{g, n}$ for all $g$
and $n$.

 Among other syzygetic ways of producing divisors on
$\mm_{g, n}$ we single the one using the Minimal Resolution
Conjecture (cf. Theorem \ref{mrc}), which can be thought of as a generalization of the divisor of higher Weierstrass points and which is especially useful in the case of
a large number of marked points. An immediate application of the
calculations in Section 4 is the following result about the Kodaira
type of $\mm_{g, n}$:

\begin{theorem}\label{mgn}
For integers $g=4, \ldots, 21$, the moduli space $\mm_{g, n}$ is of
general type for all $n\geq f(g)$ where $f(g)$ is described in the
following table.
\begin{center}
\begin{tabular}{c|cccccccccccccccccc}
$g$ & 4& 5& 6& 7& 8& 9& 10& 11& 12& 13& 14& 15& 16& 17& 18& 19& 20&
21\\
\hline $f(g)$ & 16 & 15&16& 15 &14 & 13& 11& 12& 13& 11& 10& 10& 9&
9&
9& 7& 6& 4\\
\texttt{}\end{tabular}
\end{center}
\end{theorem}
This result represents an improvement of Logan's Theorem 5.1  for
$g=4-6, 10,$ $14-16, 18-22$, the entries for the remaining values of
$g$ being those from \cite{Log}.

\noindent \textbf{Acknowledgment:} I am grateful to Sean Keel for
many discussions over the years on topics related to this circle of
ideas.

\section{Constructing divisors of small slope using syzygies}

For a projective variety $X$ and a line bundle $L$ on $X$ we denote
by $K_{i, j}(X, L)$ the Koszul cohomology group obtained from the
complex
$$\wedge^{i+1} H^0(L)\otimes H^0(L^{\otimes
(j-1)})\longrightarrow \wedge^i H^0(L)\otimes H^0(L^{\otimes
j})\longrightarrow \wedge^{i-1} H^0(L)\otimes H^0(L^{\otimes
(j+1)}),$$ where the maps are the Koszul differentials (cf.
\cite{Gr}). Assume $L$ is globally generated and $M_L$ is
the vector bundle on $X$ defined by the exact sequence
$$ 0\rightarrow M_L\rightarrow H^0(L)\otimes \OO_X\rightarrow L\rightarrow
0.$$  A simple argument using the exact sequences
$$ 0\longrightarrow \wedge^a M_L\otimes L^{\otimes b} \rightarrow
\wedge^a H^0(L)\otimes L^{\otimes b}\longrightarrow \wedge^{a-1}
M_L\otimes L^{\otimes (b+1)}\longrightarrow 0$$ for various $a$ and
$b$, shows that there is an identification
\begin{equation}\label{koszul}
 K_{a, b}(X, L)=\frac{H^0(\wedge^a
M_L\otimes L^{\otimes b})}{\mbox{Image} \{\wedge^{a+1} H^0(L)\otimes
H^0(L^{\otimes (b-1)})\}}\ .
\end{equation}
\noindent From now on we fix integers $i\geq 0$ and $s\geq 1$ and set
$$r:=2s+si+i, g:=rs+s, \mbox{ and } d:=rs+r.$$ We introduce the open substack  $\cM_g^0$
 of $\cM_g$ corresponding to curves $[C]\in \cM_g$ such that
$W^r_{d-1}(C)=\emptyset$ and $W^{r+1}_d(C)=\emptyset$. Then
$\mbox{codim}(\cM_g-\cM_g^0, \cM_g)\geq 2$. We denote by
$\mathfrak{Pic}^d$ the degree $d$ Picard stack over $\cM_g$
(precisely, the \'etale sheafification of the Picard functor). In
particular if $\rm{Pic}^d_{M_g}$ is the coarse moduli space
associated to $\mathfrak{Pic}^d$, then for any $\cM_g$-scheme
$T\rightarrow \cM_g$ originating from a family of genus $g$
 curves $\mathcal{X}\rightarrow T$, the fibre product $T\times_{M_g} \rm{Pic}^d_{M_g}$ is the relative Picard
 algebraic space $\rm{Pic}^d_{\mathcal{X}/T}$. We denote by
$\mathfrak G^r_d\subset \mathfrak{Pic}^d$ the stack parameterizing
pairs $[C, L]$ with $[C]\in \cM_g$ and $L\in W^r_d(C)$ and by
$\sigma:\mathfrak G^r_d\rightarrow \cM_g$ the natural projection.
Since $\rho(g, r, d)=0$,  by general Brill-Noether theory, the
general curve of genus $g$ will have finitely many $\mathfrak
g^r_d$'s and there exists a unique irreducible component of
$\mathfrak G^r_d$ which maps onto $\cM_g$. Moreover, the image of
any component
 of $\mathfrak{G}^r_d$ having dimension $\geq 3g-2$ is a substack of codimension $\geq 2$ in $\cM_g$ (cf.
 Corollary \ref{onecomp}), thus one can ignore these extraneous components of $\mathfrak{G}^r_d$ when doing
 divisor class calculations on $\mm_g$.

We shall define a determinantal substack of $\mathfrak G^r_d$
consisting of those pairs $[C, L]$ satisfying the condition $K_{i,
2}(C, L)\neq 0$. We denote by $\pi:\cM_{g, 1}^0\rightarrow \cM_g^0$ the
universal curve and by $\L$ a universal Poincar\'e bundle on the
fibre product $\cM_{g, 1}^0\times_{\cM_g^0}\mathfrak G^r_d$ (In the case such an$\L$ does not exist, we pass to an \'etale surjection $\Sigma\rightarrow \mathfrak{G}^r_d$ such that $\Sigma$ is a scheme and $\cM_{g, 1}^0\times_{\cM_g^0} \Sigma$ admits
a Poincar\'e bundle and we carry out the construction at this level. In the end
our construction does not depend on the choice of $\Sigma$, see also \cite{Est2}, Section 6.2).   If
$p_1:\cM_{g, 1}^0\times_{\cM_g^0} \mathfrak G^r_d \rightarrow \cM_{g,
1}^0$ and $p_2:\cM_{g,1}^0\times_{\cM_g^0} \mathfrak G^r_d\rightarrow
\mathfrak G^r_d$ are the natural projections, then
$\mathcal{E}:=p_{2 *}( \L)$ is a vector bundle of rank $r+1$
and there is a tautological embedding of the pullback of the
universal curve $\cM_{g, 1}^0\times_{\cM_g^0}\mathfrak G^r_d$ into the
projective bundle $u:\PP(\E)\rightarrow \mathfrak G^r_d$. We define
the vector bundle $\mathcal{F}$  on $\PP(\E)$ by the sequence
$$0\longrightarrow \mathcal{F}\longrightarrow u^*(\E)\longrightarrow \OO_{\PP(\E)}(1)\longrightarrow 0,$$
and we further introduce two vector bundles $\cA$ and $\cB$ over
$\mathfrak G^r_d$ by setting
$$\cA:=u_*\Bigl(\wedge^i \mathcal{F}\otimes \OO_{\PP(\E)}(2)\Bigr), \mbox{ and }
\cB:=u_*\Bigl(\wedge ^i \mathcal{F}\otimes \OO_{\cM_{g, 1}^0\times
_{\cM_g^0} \mathfrak G^r_d}(2)\Bigr).$$ If
$C\stackrel{|L|}\rightarrow  \PP^r$ is the map
corresponding to a point $[C, L]\in \mathfrak G^r_d$, then
$$\cA(C,L)=H^0(\PP^r, \wedge^i M_{\PP^r}(2)) \mbox{
}\mbox{ and }\cB(C,L)=H^0(C, \wedge^i M_L\otimes L^2)$$ and there is
a vector bundle morphism $\phi:\cA\rightarrow \cB$ given by
restriction. Grauert's Theorem guarantees that both $\cA$ and $\cB$
are vector bundles over $\mathfrak G^r_d$  and their ranks are
$$\mbox{rank}(\cA)=(i+1){r+2 \choose i+2}\mbox{  and  }\mbox{
rank}(\cB)={r\choose i}\Bigl(-\frac{id}{r}+2d+1-g\Bigr)$$ (We use
that $M_L$ is a stable vector bundle, see \cite{F3},  Proposition
2.1 and this  implies that $H^1(\wedge^i M_L\otimes L^{\otimes
2})=0$, hence $\mbox{rank}(\cB)$ can be computed from Riemann-Roch).
Because of the way we chose $g, r$ and $d$ we can see that
$\mbox{rank}(\cA)=\mbox{rank}(\cB)$.

 While the construction of $\cA$
and $\cB$ clearly depends on the choice of the Poincar\'e bundle
$\mathcal{L}$ (and of $\Sigma$), it is easy to check that the vector
bundle $Hom_{\OO_{\mathfrak G^r_d}}(\cA, \cB)$ on $\mathfrak G^r_d$ as well as the morphism
$\phi \in H^0\bigl(\mathfrak G^r_d, Hom _{\OO_{\mathfrak G^r_d}}(\cA, \cB)\bigr)$ are independent of
such choices. More precisely, let us denote by $\Xi$ the collection
of pairs $\alpha:=(\pi_{\alpha}, \mathcal{L}_{\alpha})$ where
$\pi_{\alpha} :\Sigma_{\alpha}\rightarrow \mathfrak G^r_d$ is an
\'etale surjective morphism from a scheme $\Sigma_{\alpha}$ and
$\mathcal{L}_{\alpha}$ is a Poincar\'e bundle on $p_{2, \alpha}:
\cM_{g, 1}^0\times_{\cM_g^0} \Sigma_{\alpha}\rightarrow
\Sigma_{\alpha}$. Recall that if $\Sigma \rightarrow \mathfrak
G^r_d$ is an \'etale surjection from a scheme and $\mathcal{L}$ and
$\mathcal{L}'$ are two Poincar\'e bundles on $p_2: \cM_{g,
1}^0\times_{\cM_g^0}\Sigma \rightarrow \Sigma$, then the sheaf
$\mathcal{N}:=p_{2 *} Hom (\mathcal{L}, \mathcal{L'})$ is invertible
and there is a canonical isomorphism $\mathcal{L}\otimes
p_2^*\mathcal{N}\cong \mathcal{L}'$. For every $\alpha\in \Xi$ we
construct the morphism between vector bundles of the same rank
$\phi_{\alpha}:\cA_{\alpha}\rightarrow \cB_{\alpha}$ over $\Sigma_{\alpha}$ as above. Then
since a straightforward cocycle condition is met, we find that there
exists a vector bundle $Hom _{\OO_{\mathfrak{G}^r_d}}(\cA, \cB)$ on $\mathfrak G^r_d$ together
with a section $\phi \in H^0\bigl(\mathfrak{G}^r_d, Hom _{\OO_{\mathfrak{G}^r_d}} (\cA, \cB)\bigr)$ such
that for every $\alpha=(\pi_{\alpha}, \mathcal{L}_{\alpha})\in \Xi$
we have that $\pi_{\alpha}^* (Hom _{\OO_{\mathfrak{G}^r_d}} (\cA, \cB))=Hom _{\OO_{\Sigma_{\alpha}}}(\cA_{\alpha},
\cB_{\alpha})$ and $\pi_{\alpha}^*(\phi)=\phi_{\alpha}$.

\begin{theorem}\label{ni}
The cycle $\mathcal{U}_{g, i}:=\{(C,L)\in \mathfrak G^r_d :K_{i,
2}(C,L) \neq 0 \}$ is the degeneracy locus of vector bundle map
$\phi:\cA\rightarrow \cB$ over $\mathfrak G^r_d$.
\end{theorem}

\begin{proof} Along the same lines as the proof of Proposition 2.5
in \cite{F3}.
\end{proof}

Thus $\mathcal{Z}_{g, i}:=\sigma_*(\mathcal{U}_{g, i})$ is a virtual
divisor on $\cM_g$ when $g=s(2s+si+i+1)$.

\begin{remark}
Using (\ref{koszul}) it is easy to prove that for every $(C, L)\in
\mathfrak{G}^r_d$ one have the vanishing of Koszul cohomology groups
$K_{a, 0}(C, L)=0$ for all $a\geq 1$ and $K_{a, b}(C, L)=0$ for all
$b\geq 3$. Thus the only non-trivial Koszul type conditions one
could impose on $\mathfrak{G}^r_d$ involve  the groups $K_{a, 1}(C,
L)$ and $K_{a, 2}(C, L)$. Because $M_L$ is a stable vector bundle on
$C$, it is straightforward to show using (\ref{koszul}) that
$$\mbox{dim }K_{i, 2}(C, L)-\mbox{dim } K_{i+1, 1}(C, L)={r \choose
i}(2d-\frac{id}{r}+1-g)-(i+1){r+2 \choose i+2}.$$ For our choices of
$g, r$ and $d$, it follows that $\mbox{dim }K_{i+1, 1}(C,
L)=\mbox{dim }K_{i, 2}(C, L)$, hence $\mathcal{U}_{g, i}$ can also
be defined as the locus where $K_{i+1, 1}(C,L)$ fails to vanish.
This shows that, at least in the case of curves, there are no other
Koszul divisors except $\mathcal{U}_{g, i}$.

\end{remark}

To prove Theorem \ref{div} we shall extend the determinantal
structure of $\mathcal{Z}_{g, i}$ over a substack of $\mm_{g}$ whose
complement has codimension $\geq 2$. We denote by
$\widetilde{\cM}_g:=\cM_g^0\cup \bigl(\cup_{j=0}^{[g/2]}
\Delta_j^0\bigr)$ the locally closed substack of $\mm_g$ obtained by
adding to $\cM_{g}^0$  the open subsets $\Delta_j^0\subset \Delta_j$
for $1\leq j\leq [g/2]$ consisting of $1$-nodal genus $g$ curves
$C\cup_y D$, with $[C]\in \cM_{g-j}$ and $[D, y]\in \cM_{j, 1}$
being Brill-Noether general curves, and the locus $\Delta_0^0\subset
\Delta_0$ containing $1$-nodal irreducible genus $g$ curves
$C'=C/q\sim y$, where $[C, q]\in \cM_{g-1}$ is a Brill-Noether
general pointed curve and $y\in C$, together with their
degenerations consisting of unions of a smooth genus $g-1$ curve and
a nodal rational curve. One can then extend the finite covering
$\sigma:\mathfrak G^{r}_{d}\rightarrow \cM_g^0$ to a proper,
generically finite map
$$\sigma: \widetilde{\mathfrak G}^{r}_{d} \rightarrow \widetilde{\cM}_g$$ by
letting $\widetilde{\mathfrak G}^{r}_{d}$ be the space of limit
$\mathfrak g^{r}_{d}$'s on the curves from $\widetilde{\cM}_g$ which
are all tree-like (see \cite{EH1}, Theorem 3.4 for the construction
of the variety of limit linear series and also \cite{Oss} for a more
functorial approach which in the case $\rho(g, r, d)=0$ leads to the
Eisenbud-Harris space). Strictly speaking, Eisenbud and Harris have
only constructed the space of \emph{refined} limit $\mathfrak
g^r_d$'s. Using the observation that when $\rho(g,r, d)=0$ every
\emph{crude} non-refined limit $\mathfrak g^r_d$ on a curve of
compact type $C\cup _y D$, where $[C]\in \cM_j$ and $[D]\in
\cM_{g-j}$  can be canonically interpreted as a refined limit
$\mathfrak g^r_d$ on the pre-stable curve $C\cup _{y_1} \PP^1
\cup_{y_2}D$ obtained from $C\cup_y D$ by inserting a single $\PP^1$
at the node $y$, their construction can be easily adapted to cover
the case of crude $\mathfrak g^r_d$'s as well.  Note that since all
limit $\mathfrak g^r_d$'s are \emph{dimensionally proper} (cf.
\cite{EH1}, Corollary 3.7), every limit linear series from
$\widetilde{\mathfrak G}^r_d$ is smoothable.

To compute the class $[\overline{\mathcal{Z}}_{g, i}]$, we  intersect $\overline{\mathcal{Z}}_{g, i}$
 with  test curves in the boundary of $\mm_{g}$ which are defined
as follows: we fix a Brill-Noether general curve $C$ of genus $g-1$,
a general point $q\in C$ and a general elliptic curve $E$. We define
two $1$-parameter families
\begin{equation}\label{testcurves}
C^0:=\{C/y\sim q: y\in C\}\subset \Delta_0 \subset \mm_{g}
\mbox{ and }C^1:=\{C\cup _y E: y\in C\}\subset \Delta_1\subset
\mm_{g}.
\end{equation}
 It is well-known that
these families intersect the generators of
$\mbox{Pic}(\mm_{g})$ as follows:
$$
C^0\cdot \lambda=0,\ C^0\cdot \delta_0=-(2g-2), \ C^0\cdot
\delta_1=1 \mbox{ and } C^0\cdot \delta_a=0\mbox{ for }a\geq 2,
\mbox{ and}$$
$$C^1\cdot \lambda=0, \ C^1\cdot \delta_0=0, \ C^1\cdot
\delta_1=-(2g-4), \ C^1\cdot \delta_a=0 \mbox{ for }a\geq 2.$$ Next,
we fix  $2\leq j\leq [g/2]$, a general curve $C$ of genus $j$ and a
general curve pointed curve $(D, y)$ of genus $g-j$. We define the
$1$-parameter family $C^j:=\{C\cup_y D: y\in C\}\subset
\Delta_j\subset \mm_{g}$. We have that
$$C^j\cdot \lambda=0, \ C^j\cdot \delta_a=0 \mbox{ for }a\neq j
\mbox{ and } C^j\cdot \delta_j=-(2j-2).$$ We  review the notation
used in the theory of limit linear series
 (see \cite{EH1}).
 If $X$ is a
tree-like curve and $l$ is a limit $\mathfrak g^r_d$ on $X$, for a
component  $Y$ of $X$ we denote by $l_Y=(L_Y, V_Y\subset H^0(L_Y))$
the $Y$-aspect of $l$. For a point $y\in Y$ we denote by
$\{a^{l_Y}_i(y)\}_{i=0, \ldots, r}$ the \emph{vanishing sequence} of
$l$ at $y$, by $\{\alpha^{l_Y}_i(y)=a_i^{l_Y}(y)-i\}_{i=0, \ldots,
r}$ the \emph{ramification sequence} and by $\rho(l_Y, y):=\rho(g,
r, d)-\sum_{i=0}^r \alpha^{l_Y}_i(y)$ the adjusted Brill-Noether
number with respect to $y$.
\begin{proposition}\label{limitlin}
(1) Let $C_y^1=C\cup_y E$ be an element of $\Delta_1^0$. If $(l_C,
l_E)$ is a limit  $\mathfrak g^{r}_{d}$ on $C_y^1$, then
$V_C=H^0(L_C)$ and $L_C\in W^{r}_{d}(C)$ has a cusp at $y$. If $y\in
C$ is a general point, then $l_E=\bigl(\OO_E(dy),
(d-r-1)y+|(r+1)y|\bigr)$, that is, $l_E$ is uniquely determined. If
$y\in C$ is one of the finitely many points for which there exists
$L_C\in W^{r}_{d}(C)$ such that $\rho(L_C, y)=-1$, then
$l_E(-(d-r-2)y)$ is a $\mathfrak g_{r+2}^{r}$ with vanishing
sequence at $y$ being $\geq (0, 2, 3, \ldots, r, r+2)$.  Moreover,
at the level of $1$-cycles we have the identification
$\sigma^*(C^1)\equiv X + \nu \ T$, where
$$X:=\{(y, L)\in C\times W^r_d(C) :h^0(C, L(-2y))\geq r\},$$
$T\cong \PP\bigl(H^0(\OO_E((r+2)y))/H^0(\OO_E(ry))\bigr)$ is the curve consisting of $\mathfrak g^{r}_{r+2}$'s on $E$
with vanishing $\geq (0, 2, \ldots, r, r+2)$ at the fixed point
$y\in E$ and $\nu$ is an explicitly known positive integer.

\noindent (2) Let $C_y^0=C/y\sim q$ be an element of $\Delta_0^0$.
Then limit linear series of type $\mathfrak g^{r}_{d}$ on $C_y^0$
are in 1:1 correspondence with complete linear series $L$ on $C$ of
type $\mathfrak g^{r}_{d}$ satisfying the condition $h^0(C, L\otimes
\OO_C(-y-q))=h^0(C,L)-1.$ There is an isomorphism between the
cycle $\sigma^*(C^0)$ of $\mathfrak g^{r}_{d}$'s on all curves
$C_y^0$ with $y\in C$ and the smooth curve
$$Y:=\{(y, L)\in C\times W^r_d(C): h^0(C, L(-y-q))\geq r\}.$$
\end{proposition}

\begin{proof}
Part (1) is  similar to the proof of Proposition 3.3 from \cite{F3} and
we omit the details. For part (2), we claim that for any limit $\mathfrak g^r_d$ on a curve $C_y^0$ where $y\in C$, the underlying torsion free sheaf is actually locally free. Indeed, otherwise the underlying sheaf would be of the form $\nu_*(L)$, where $\nu:C\rightarrow C^y_0$ is the normalization map and
$L\in W^r_{d-1}(C)$. But $[C]\in \cM_{g-1}$ is assumed to be Brill-Noether general, hence $W^r_{d-1}(C)=\emptyset$.
\end{proof}
Throughout this paper we routinely use basic facts from Schubert
calculus which we briefly recall. If $\GG(r, d)$ denotes the
Grassmannian of $r$-planes in $\PP^d$ and $$\mathbb
C^{d+1}=V_0\supset V_1 \supset\cdots \supset V_{r+1}=0$$ is a
decreasing flag, then for any Schubert index $0\leq \alpha_0\leq
\ldots \leq \alpha_{r}\leq d-r$, we define the Schubert cycle
$$\sigma_{(\alpha_0, \ldots, \alpha_r)}:=\{\Lambda \in \GG(r, d):
\mbox{dim}(\Lambda\cap V_{\alpha_i+i})\geq r+1-i, \mbox{ for } i=0,
\ldots, r\}.$$ (This differs slightly from the standard notation
from e.g. \cite{FuPr}, but it seems better suited for dealing with
ramification sequences of linear series). Often we use the fact that
if $(\alpha_0, \ldots, \alpha_r)$ is a Schubert index and $g$ is an
integer such that $rg+\sum_{i=0}^r \alpha_i=(r+1)(d-r)$, then there
is an identity in $H^*(\GG(r, d))$:
\begin{equation}\label{schubert0}
\sigma_{(\alpha_0, \ldots, \alpha_r)}\cdot \sigma_{(0, 1, \ldots,
1)}^g=g!\frac{\prod_{i<j} (\alpha_j-\alpha_i+j-i)}{\prod_{i=0}^r
(g-d+i+\alpha_i+r)!}.
\end{equation}

\begin{proposition}\label{limitlinj}
 Let $[C]\in \cM_j$ be a general curve with $g-2\geq j\geq [g/2]$
and $C^j\subset \Delta_j\subset \mm_g$ the associated test curve of type
$(j, g-j)$. Then one has the following equality of $1$-cycles in $\widetilde{\mathfrak G}^r_d$:
$$\sigma^*(C^j)=\sum_{(\alpha_0, \ldots, \alpha_r)\in \mathcal{P}_1} N_{g-j, \alpha} \cdot X_{j, \alpha}
+\sum_{(\beta_0, \ldots, \beta_r)\in \mathcal{P}_2} M_{j,
\beta}\cdot Y_{g-j, \beta}+\sum_{(\beta_0, \ldots, \beta_r)\in
\mathcal{P}_3} Q_{g-j, \beta}\cdot U_{j, \beta},
$$
where  we introduce the following notations:
$\mathcal{P}_1:=\{(0\leq \alpha_0\leq \ldots \leq \alpha_r\leq s):
\sum_{i=0}^r \alpha_i=j\},$
\begin{center} $\mathcal{P}_2:=\{(0\leq
\beta_0\leq \ldots \leq \beta_r\leq s+1):\sum_{i=0}^r \beta_i=j+1,
\beta_{r-1}\leq s\}$,\end{center}
\begin{center}
$\mathcal{P}_3:=\{(0=\beta_0<\beta_1\leq \ldots \leq \beta_r\leq
s+1):\sum_{i=0}^r \beta_i=r+1+j\}$,
\end{center}
 \begin{center}
 $M_{j, \beta}:=\sigma_{(0, 1, \ldots, 1)}^j\cdot \sigma_{(\beta_0, \ldots, \beta_r)}\in H^*(\GG(r,
r+j)) \mbox{ for }\beta\in \mathcal{P}_2,$
\end{center}
\begin{center} $N_{g-j, \alpha}:=\sigma_{(0, 1,
\ldots, 1)}^{g-j}\cdot \sigma_{(j-\alpha_r, \ldots, j-\alpha_0)} \in
 H^*(\GG(r, d)) \mbox{ for } \alpha\in \mathcal{P}_1,$
\end{center}
\begin{center}
$Q_{g-j, \beta}:=\sigma_{(0, 1, \ldots, 1)}^{g-j}\cdot
\sigma_{(j+1-\beta_r,
  \ldots, j+1-\beta_1, j+1)}\in H^*(\GG(r, d))$ for $\beta\in
\mathcal{P}_3,$
  \end{center}
\begin{center}
$X_{j, \alpha}:=\{(y, L_C)\in C\times \rm{Pic}$$^{r+j}(C):
\alpha^{L_C}_i(y)\geq \alpha_i \mbox{ for } i=0\ldots r\},\mbox{
}\alpha \in \mathcal{P}_1$,
\end{center}
\begin{center}
 $Y_{g-j, \beta}:=\{l_D\in G^r_d(D): \alpha^{l_D}_i(y)\geq
j-\beta_{r-i}\mbox{ for }i=0\ldots r\}, \mbox{ }\beta\in
\mathcal{P}_2$,
\end{center}
\begin{center}
$U_{j, \beta}:=\{(y, l_C)\in C\times G^r_{r+j+1}(C):
\alpha^{l_C}(y)\geq (0, \beta_1, \ldots, \beta_r)\}, \mbox{ for
}\beta \in \mathcal{P}_3$.
\end{center}
\end{proposition}
\begin{proof} Suppose that $l=(l_C, l_D)$ is a limit $\mathfrak g^r_d$ on $C\cup _y D$. It is easy
to see that the generic point of any component of $\sigma^*(C^j)$ corresponds to a refined
limit $\mathfrak g^r_d$,
so we may assume that
$l$ is refined as well. If $(\alpha_0, \ldots, \alpha_r)$ is the ramification sequence of $l_C$ at $y$,
 then the condition that  $[D, y]\in \cM_{g-j, 1}$ carries a $\mathfrak g^r_d$
with ramification sequence at $y$ being at least $(d-r-\alpha_r,
\ldots, d-r-\alpha_0)$, is that $\sigma_{(0, 1, \ldots,
1)}^{g-j}\cdot \sigma_{(d-r-\alpha_r, \ldots, d-r-\alpha_0)}\neq 0
\in H^*(\GG(r, d))$. Using the Littlewood-Richardson rule, we find
that this implies that $\alpha_r\leq rs+s-j$. A similar reasoning
can be used for $C$. Degenerating $C$ to a stable curve consisting
of a rational spine and $j$ elliptic tails, we obtain that if there
exists a point $y\in C$ and a $\mathfrak g^r_d$ with ramification
sequence $(\alpha_0, \ldots, \alpha_r)$ at $y$, then either $y$
specializes to a point on the rational spine in which case we find
the condition $\sigma_{(0, 1, \ldots, 1)}^j\cdot \sigma_{(\alpha_0,
\ldots, \alpha_r)}\neq 0 \in H^*(\GG(r, d))$ which implies that
$\alpha_0\geq rs-j$, or else, $y$ specializes to a point on one of
the elliptic tails in which case we find that there must exist two
indices $0\leq e<f\leq r$ with $\alpha_e\geq \alpha_{e-1}+1$ and
$\alpha_{f}\geq \alpha_{f-1}+1$, such that $\sigma_{(0, 1, \ldots,
1)}^{j-1}\cdot \sigma_{(\alpha_0+1, \alpha_{e-1}+1, \alpha_e,
\alpha_{e+1}+1, \ldots, \alpha_{f-1}+1, \alpha_f, \alpha_{f+1}+1,
\ldots, \alpha_r+1)}\neq 0$. This last condition leads to the
inequality $\alpha_0\geq \mbox{max}\{0, rs-j-1\}$.

Suppose we are in the first case, that is, $\alpha_0\geq rs-j$ and
moreover $\rho(l_C, y)=\rho(l_D, y)=0$, which is the situation which
occurs for a generic choice of $y\in C$. Then $l_C(-(rs-j)y)=|L_C|$,
where $L_C\in \mbox{Pic}^{r+j}(C)$ with $\alpha_r^{L_C}(y)\leq s$
and $\sum_{i=0}^r \alpha_i^{L_C}(y)=j$, that is, $(L_C, y)\in X_{j,
\alpha}$. If $\alpha_0\geq rs-j$ but now $\rho(l_C, y)=-1$ and
$\rho(l_D, y)\leq 1$, then $\{\alpha_i^{L_C}(y)-(rs-j)\}_{i=0\ldots
r}$ must be one of the partitions from the set $\mathcal{P}_2$.
Choosing such a partition, we have $M_{j, \beta}$ choices for the
$C$-aspect, while $l_D\in Y_{g-j, \alpha^{L_C}(y)-(rs-j)}$. Finally
let us assume that we are in the case $\alpha_0^{l_C}(y)=rs-j-1$.
Then necessarily $\alpha_1^{l_C}(y)\geq rs-j$, $\rho(l_C,
y)=\rho(l_D, y)=0$ and $l_C(-(rs-j-1)y) \in U_{j, \beta}$, where
$\beta_i:=\alpha_i^{l_C}(y)-(rs-j-1)$ for $i=0\ldots r$. This
accounts for the third sum in $\sigma^*(C^j)$. Arguing along the
lines of \cite{EH2}, Lemma 3.4, $\widetilde{\mathfrak G}^r_d$ is
smooth along $\sigma^*(C^j)$ and since all limit $\mathfrak g^r_d$
described in this proof are smoothable, we obtain that the claimed
formula holds at the level of $1$-cycles (including multiplicities).
\end{proof}
The next corollary shows that "ghost" components of
$\widetilde{\mathfrak G}^r_d$ having dimension $> 3g-3$, do not
matter in the calculation of $[\sigma_*(\G_{i, 2}-\H_{i, 2})]$.
\begin{corollary}\label{onecomp}
In the case $\rho(g, r, d)=0$, every irreducible component of
$\mathcal{Z}$ of $\widetilde{\mathfrak G}^r_d$ such that
$\rm{dim}$$(\mathcal{Z})\geq 3g-2$ has the property that $\rm{dim
}$$\ \sigma(\mathcal{Z})\leq 3g-5$.
\end{corollary}
\begin{proof} For a general $[C]\in \cM_g$, the scheme $W^r_d(C)$ is reduced and $0$-dimensional, thus every component of $\widetilde{\mathfrak G}^r_d$ mapping dominantly onto $\widetilde{\cM}_g$ must have dimension $3g-3$. Suppose that $\mathcal{Z}$ is a component of dimension at least $3g-2$ such that $\sigma(\mathcal{Z})$ is a divisor on $\widetilde{\cM}_g$. Then for any $[C]\in \sigma(\mathcal{Z})$ we have that $\mbox{dim}\ \sigma^{-1}([C])\geq 2$. Since this property does not hold along any of the curves $\sigma^*(C^j)$ for any $0\leq j\leq [g/2]$ (see Proposition \ref{limitlinj}), it follows that
$\sigma(\mathcal{Z})$ is disjoint from the test curves $C^j\subset \mm_g$ for all $j\geq 0$. This implies then that $[\sigma(\mathcal{Z})]=0\in \mbox{Pic}(\widetilde{\cM}_g)$, hence $\sigma(\mathcal{Z})=0$ (use that the Satake compactification of $\cM_g$ has boundary of codimension $2$). This is a contradiction.
\end{proof}

Let $C$ be a Brill-Noether general curve of genus $g-1$ (recall that
$g=rs+s$ and $d=rs+r$). Then $\mbox{dim } W^r_d(C)=r$ and it is easy
to see that $C$ carries no $\mathfrak g^r_{d-1}$'s or $\mathfrak
g^{r+1}_d$'s, hence every $L\in W^r_d(C)$ corresponds to a complete
and base point free linear series.  We denote by $\L$ a Poincar\'e
bundle on $C\times \mbox{Pic}^d(C)$ and by $\pi_1:C\times
\mbox{Pic}^d(C)\rightarrow C$ and $\pi_2:C\times
\mbox{Pic}^d(C)\rightarrow \mbox{Pic}^d(C)$ the projections. We
define the cohomology class $\eta=\pi_1^*([\mbox{point}])\in
H^2(C\times \mbox{Pic}^d(C))$, and if $\delta_1,\ldots,
\delta_{2g}\in H^1(C, \mathbb Z)\cong H^1(\mbox{Pic}^d(C), \mathbb
Z)$ is a symplectic basis, then we set
$$\gamma:=-\sum_{\alpha=1}^g
\Bigl(\pi_1^*(\delta_{\alpha})\pi_2^*(\delta_{g+\alpha})-\pi_1^*(\delta_{g+\alpha})\pi_2^*(\delta_
{\alpha})\Bigr).$$ We have the formula (cf. \cite{ACGH}, p. 335)
$c_1(\L)=d\eta+\gamma,$ corresponding to the Hodge decomposition of
$c_1(\L)$. We also record that $\gamma^3=\gamma \eta=0$, $\eta^2=0$
and $\gamma^2=-2\eta \pi_2^*(\theta)$. On $W^r_d(C)$  we have the
tautological rank $r+1$ vector bundle
$\mathcal{E}:=(\pi_2)_{*}(\mathcal{L}_{| C\times W^r_d(C)})$. The
Chern numbers of $\mathcal{E}$ can be computed using the Harris-Tu
formula (cf. \cite{HT}): if we write $\sum_{i=0}^r
c_i(\mathcal{E}^{\vee})=(1+x_1)\cdots (1+x_{r+1})$, then for every
class $\zeta \in H^*(\mbox{Pic}^d(C), \mathbb Z)$ one has the formula
\footnote{There is a confusing sign error in the formula (1.4) in
\cite{HT}: the formula is correct as it is appears in \cite{HT}, if
the $x_j$'s denote the Chern roots of the \emph{dual} of the kernel
bundle.}
$$x_1^{i_1}\cdots x_{r+1}^{i_{r+1}}\
\zeta=\mbox{det}\Bigl(\frac{\theta^{g-1+r-d+i_j-j+l}}{(g-1+r-d+i_j-j+l)!}\Bigr)_{1\leq
j, l\leq r+1}\ \zeta.$$ If  we use the expression of the Vandermonde
determinant, we get the formula
$$\mbox{det}\Bigl(\frac{1}{(a_j+l-1)!}\Bigr)_{1\leq j, l\leq
r+1}=\frac{\Pi_{ j>l}\ (a_l-a_j)}{\Pi_{j=1}^{r+1}\ (a_j+r)!}.$$ By
repeatedly applying this formula we compute all the intersection
numbers on $W^r_d(C)$ which we shall need:
\begin{lemma}\label{vandermonde}
If $c_i:=c_i(\mathcal{E}^{\vee})$ we have the following identities
in $H^*(W^r_d(C), \mathbb Z)$:

$ (1) \ \ \mbox{ } c_{r-1}\theta=\frac{r(s+1)}{2} c_r$

$ (2)  \ \ \mbox{ } c_{r-2} \theta^2=\frac{r(r-1)(s+1)(s+2)}{6} c_r$

$ (3)  \ \ \mbox{ } c_{r-2} c_1
\theta=\frac{r(s+1)}{2}\bigl(1+\frac{(r-2)(r+2)(s+2)}{3(s+r+1)}\bigr)
c_r$

$ (4) \ \ \mbox{ } c_{r-1}
c_1=(1+\frac{(r-1)(r+2)(s+1)}{2(s+r+1)})c_r$

$ (5) \ \ \mbox{ } c_r=\frac{1! \ 2!\cdots (r-1)!\ (r+1)!}{(s-1)!\
(s+1)!\ (s+2)! \cdots (s+r)!}\theta^{g-1}\ .$
\end{lemma}

We point out that the constant $c_r$ equals the number of linear
series $\mathfrak g^r_d$ on a general curve of genus $g$ (note that
$\rho(g, r, d)=0$). In Section 3 we shall use the following result:
\begin{lemma}\label{h1} If $[C]\in \cM_{g-1}$, then one has the following identity in $H^*(W^r_d(C), \mathbb Z)$:
$$c_1\bigl(R^1 \pi_{2  *}(\mathcal{L}_{| C\times
W^r_d(C)})\bigr)=\theta-c_1(\mathcal{E}^{\vee}).$$
\end{lemma}
\begin{proof} Let us recall how one can obtain a determinantal structure on
$W^r_d(C)$. Once we fix a divisor $D\in C_e$ of degree $e>>0$,
$W^r_d(C)$ is the degeneracy locus of rank $d+e-g+1-r$ of the vector
bundle map $(\pi_2)_*\bigl(\mathcal{L}\otimes \OO(\pi_1^*D)\bigr)
\rightarrow (\pi_2)_*\bigl(\mathcal{L}\otimes \OO(\pi_1^*D)_{|
\pi_1^*(D)}\bigr)$. Consequently, we have an exact sequence of
vector bundles over $W^r_d(C)$:
$$0\longrightarrow \mathcal{E}\longrightarrow
(\pi_{2})_*\bigl(\mathcal{L}\otimes
\OO(\pi_1^*D)\bigr)\longrightarrow
(\pi_{2})_*\bigl(\mathcal{L}\otimes \OO(\pi_1^*D)_{| \pi_1^*
D}\bigr)\longrightarrow R^1\pi_{2
*}(\mathcal{L}_{| C\times W^r_d(C)})\longrightarrow 0,$$ from which the claim follows by
using that $(\pi_2)_*\bigl(\mathcal{L}\otimes \OO(\pi^*D)_{|
\pi_1^*D}\bigr)$ is numerically trivial while
$c_t\bigl((\pi_2)_*(\mathcal{L}\otimes
\OO(\pi_1^*D))\bigr)=e^{-\theta}$ (cf. \cite{EH3} or \cite{ACGH}).
\end{proof}

For integers $0\leq a\leq r$ and $ b\geq 2$ we shall define vector
bundles $\G_{a,b}$ and $\H_{a, b}$ over $\sigma^{-1}(\cM_g^0 \cup
\Delta_{0}^0\cup \Delta_1^0)\subset \widetilde{\mathfrak G}^{r}_{d}$
which over the locus corresponding to smooth curves have fibres
$$\G_{a,b}(C,L)=H^0(C, \wedge^a M_L\otimes L^{\otimes b})\ \mbox{ and }\ \H_{a, b}(C,L)=H^0(\PP^r,
\wedge^a M_{\PP^r}(b))$$ for each $(C, L)\in \mathfrak G^{r}_{d}$
giving a map $C\stackrel{|L|} \rightarrow \PP^r$.  Clearly $\G_{i, 2
| \mathfrak G^{r}_{d}}=\mathcal{B}$ and $\H_{i, 2 | \mathfrak
G^r_d}=\cA$, where $\cA$ and $\mathcal{B}$ are the vector bundles
introduced in Theorem \ref{ni}. Partially extending these bundles
over the boundary of $\widetilde{\mathfrak G}^r_d$ will enable us to
compute the $\lambda, \delta_0$ and $\delta_1$ coefficients of
$\overline{\mathcal{Z}}_{g,i}$ and determine the slope $s(\overline{\mathcal{Z}}_{g, i})$.

\begin{proposition} For each $b\geq 2$ there
exists a vector bundle $\G_{0,b}$ over $\sigma^{-1}(\cM_{g}^0\cup
\Delta_0^0\cup \Delta_1^0) \subset \widetilde{\mathfrak G}^{r}_{d}$
having rank $bd+1-g$ whose fibres admit the following description:
\begin{itemize}
\item For $(C, L)\in \mathfrak G^{r}_{d}$, we have that
$\G_{0,b}\bigl(C,L)=H^0(C, L^{\otimes b})$.
\item For $t=(C\cup_y
E, L)\in \sigma^{-1}(\Delta_1^0)$, where $L\in W^r_d(C)$ has a cusp
at $y\in C$, we have that
$$\G_{0,b}(t)=H^0(C, L^{\otimes b}(-2y))+\mathbb C\cdot
u^b\subset H^0\bigr(C, L^{\otimes b}),$$ where $u\in H^0(C, L)$ is
any section such that $\rm{ord}$$_y(u)=0$. \item For $t=(C/y\sim q,
L)\in \sigma^{-1}(\Delta_0^0)$, where $q,y\in C$ and $L\in W^r_d(C)$
is such that $h^0(C, L(-y-q))=h^0(L)-1$, we have that
$$\G_{0,b}(t)=H^0(C, L^{\otimes b}(-y-q))\oplus \mathbb C\cdot u^b\subset H^0(C, L^{\otimes b}),$$
where $u\in H^0(C, L)$ is a section such that $\rm{ord}$$_y(u)=\rm{ord}$$_q(u)=0$.
\end{itemize}
\end{proposition}
\begin{proof} Very similar to Proposition 3.9 in \cite{F3}.
\end{proof}

Having defined the vector bundles $\G_{0, b}$ we now define
inductively all vector bundles $\G_{a, b}$ by the exact sequence
\begin{equation}\label{gi}
0\longrightarrow \G_{a, b}\longrightarrow \wedge^a \G_{0,
1}\otimes \G_{0, b}\stackrel{d_{a, b}}\longrightarrow \G_{a-1,
b+1}\longrightarrow 0.
\end{equation}
To define $\H_{a, b}$ is even easier. We set $\H_{0,
b}:=\mbox{Sym}^b \G_{0, 1}$ for all $b\geq 1$ and we define $\H_{a,
b}$ inductively via the exact sequence
\begin{equation}\label{hi}
0\longrightarrow \H_{a, b}\longrightarrow \wedge^a \H_{0, 1}\otimes
\mbox{Sym}^b \H_{0, 1}\longrightarrow \H_{a-1, b+1}\longrightarrow
0.
\end{equation}
The surjectivity of the right map in (\ref{hi}) is obvious, whereas
to prove that $d_{a, b}$ is surjective, one argues like in
\cite{F3}, Proposition 3.10. There is a natural vector bundle
morphism $\phi_{a, b}:\H_{a, b}\rightarrow \G_{a, b}$. Moreover
$\mbox{rank}(\H_{i, 2})=\mbox{rank}(\G_{i, 2})$ and the degeneracy
locus of $\phi_{i, 2}$ is the codimension one compactification of
$\mathcal{Z}_{g, i}$ over $\cM_g^0\cup \Delta_0^0\cup \Delta_1^0$.

We prove a technical result we shall use later for extending the
bundles $\G_{0, b}$ with $b\geq 2$ over the boundary of $\mathfrak G^r_d$. It can be interpreted as saying that on
a suitably general curve, the ramification points of a linear series are distinct from those of its higher order powers.
\begin{proposition}\label{vanish2}
Fix integers $s\geq 2, r\geq 2s$ and a partition $0\leq \beta_0\leq
\beta_1\leq \cdots \leq \beta_r\leq s$ such that $\sum_{i=0}^r
\beta_i=\gamma$. Let $(D, y)$ be a general pointed curve of genus $\gamma\geq 3s$. Then for every line bundle $L_D\in
\mathrm{Pic}^{\gamma+r}(D)$ satisfying the conditions
$\alpha_i^{L_D}(y)=\beta_i$ for $0\leq i\leq r$, we have that
$$H^0\bigl(D, K_D\otimes L_D^{\otimes (-2)}\otimes
\OO_D(ay)\bigr)=0, \mbox{ for all } a\leq 2(r+s).$$
\end{proposition}

\begin{proof} Clearly it suffices to prove the theorem in the case
$a=2(r+s)$. We degenerate $(D, y)$ to a stable curve $E_0\cup \ldots
\cup E_{\gamma-1}$, consisting of a string of elliptic curves such
that $E_{i-1}\cap E_i=\{p_i\}$ for $1\leq i\leq \gamma-1$. Moreover,
we assume that $y=p_0$ specializes to a point lying on $E_0$ and
that the differences $p_i-p_{i-1}\in \mbox{Pic}^0(E_{i-1})$ are not
torsion for all $1\leq i\leq \gamma-1$. We assume by contradiction that $H^0\bigl(D, K_D\otimes
L_D^{\otimes (-2)}\otimes \OO_D(2(r+s)y)\bigr)\neq 0$ for some $L_D\in \mbox{Pic}^{\gamma+r}(D)$ and denote by
$L_{E_i}\in \mbox{Pic}^{\gamma+r}(E_i)$ the $E_i$-aspect of
the induced limit linear series $\mathfrak g^r_{r+\gamma}$ on $\cup_{i=0}^{\gamma-1} E_i$ satisfying the
ramification conditions $\alpha^{L_{E_0}}_t(p_0)=\alpha_t$ for $0\leq
t\leq r$. Fix an integer $1\leq i\leq \gamma-1$. By the additivity of the Brill-Noether number, we have that
 $\rho(L_{E_{i-1}}, p_i, p_{i-1})=0$  and there exists an integer $0\leq k\leq r$ such that $\alpha_t^{L_{E_i}}(p_i)=\alpha_t^{L_{E_{i-1}}}(p_{i-1})+1$ for $t\neq k$ while $
\alpha_k^{L_{E_i}}(p_i)=\alpha_k^{L_{E_{i-1}}}(p_{i-1})$. In particular,
 $$L_{E_{i-1}}=\OO_{E_{i-1}}\bigl((\alpha_k^{L_{E_{i-1}}}(p_{i-1})+k)\cdot p_{i-1}+(\gamma+r-k-\alpha_k^{L_{E_{i-1}}}(p_{i-1}))\cdot p_i\bigr),$$ that is $L_{E_{i-1}}$ corresponds to a divisor supported only at the points $p_{i-1}$ and $p_i$. Our assumption implies that  for all $0\leq i\leq \gamma-1$ there exist
sections
$$0\neq \rho_i\in H^0\bigl(E_i,
\OO_{E_i}(2(r+s+i)p_i+2(\gamma-i-1)p_{i+1})\otimes L_{E_i}^{\otimes
(-2)})\bigr)$$ satisfying the compatibility conditions $$\mbox{ord}_{p_i}(\rho_i)\geq
\mbox{ord}_{p_{i-1}}(\rho_{i-1}) \mbox{ and
}\mbox{ord}_{p_{i}}(\rho_{i-1})+\mbox{ord}_{p_i}(\rho_{i})=2s-2.$$
We reach a contradiction once we show that
$$\mbox{ord}_{p_{\gamma-1}}(\rho_{\gamma-1})> 2s-2=\mbox{deg}(\OO_{E_{\gamma-1}}(2(r+s+\gamma-1)p_{\gamma-1}\otimes L_{E_{\gamma-1}}^{\otimes (-2)}))$$ which gives that $\rho_{\gamma-1}=0$. Assume now that $\mbox{ord}_{p_{i+1}}(\rho_{i+1})=\mbox{ord}_{p_i}(\rho_i)$ for some $0\leq i\leq \gamma-2$. Then $\mbox{ord}_{p_i}(\rho_i)+\mbox{ord}_{p_{i+1}}(\rho_i)=2s-2$, hence the section $\rho_i$ vanishes only at $p_i$ and $p_{i+1}\in E_i$ and $\mbox{ord}_{p_i}(\rho_i)=2b$ for some integer $b\geq 0$. We must have that $L_{E_i}=\OO_{E_i}((r+s+i-b)\cdot p_i+(\gamma-i-s-b)\cdot p_{i+1})$ (we use that $p_{i+1}-p_i\in \mbox{Pic}^0(E_i)$ is not torsion). In particular $r+s+i-b$ is one entry in the vanishing sequence $a^{L_{E_i}}(p_i)$ and the vanishing sequence $a^{L_{E_{i+1}}}(p_{i+1})$ is obtained from $a^{L_i}(p_i)$ by raising all entries by $1$, except for $r+s+i-b$ which remains unmodified. Obviously then, the number $r+s+i+1-b$ cannot appear in the vanishing sequence $a^{L_{E_{i+1}}
 }(p_{i+1})$. But this implies that $\mbox{ord}_{p_{i+2}}(\rho_{i+2})\geq \mbox{ord}_{p_{i+1}}(\rho_{i+1})+1$. This argument shows that as we trace the non-decreasing sequence of vanishing orders $\{\mbox{ord}_{p_i}(\rho_i)\}_{i=0}^{\gamma-1}$  along any group of $3$ consecutive components $E_{i-1}, E_i$ and $E_{i+1}$, we will find at least $2$ along which $\mbox{ord}_{p_i}(\rho_i)$ jumps. Since $\gamma\geq 3s$, we find that $\mbox{ord}_{p_{\gamma-1}}(\rho_{\gamma-1})>2s-2$ and this brings about a contradiction.
\end{proof}

Next we extend  $\G_{a, b}$ and $\H_{a, b}$ over the divisors
$\sigma^{-1}(\Delta_j^0)$ for $[g/2]\leq j\leq g-2$.
\begin{proposition}\label{ext0}
\begin{enumerate}
\item \noindent For $g=rs+s, d=rs+r$ and $b\geq 1$, there exists a vector bundle $\G_{0, b}$ defined over
$\widetilde{\mathfrak G}^r_d$, extending the already constructed
vector bundle $\G_{0, b}$ over $\sigma^{-1}(\cM_g^0 \cup \Delta_0^0
\cup \Delta_1^0)$ and such that if $t=(C\cup_y D, l_C, l_D)\in
\sigma^{-1}(\Delta_j^0)$, where $ g(C)=j\geq [s(r+1)/2]$,
$g(D)=g-j\geq 2$ and $(l_C, l_D)$ is a limit $\mathfrak g^r_d$ on
$C\cup_y D$, then
$$\G_{0, b}(t)=H^0(C\cup_y D, L^{\otimes b}),$$ with
$$L:=\Bigl(L_C=l_C(-(rs-j)y), L_D=l_D(-(j+r)y)\Bigr)\in
\mathrm{Pic}^{j+r}(C)\times \mathrm{Pic}^{rs-j}(D),$$ in the case
$s(r+1)/2 \leq j \leq s(r-1)=g-2s$, and
$$L:=\Bigl(L_C=l_C(-(g-j+r)y), L_D=l_D(-(j-s)y)\Bigr) \in
\mathrm{Pic}^{j-s}(C)\times \mathrm{Pic}^{rs+s-j+r}(D),$$ in the case
$s(r-1)<j\leq s(r+1)-2$.
\item For each $0\leq a\leq r, b\geq 1$ there exists a vector bundle $\H_{a, b}$
over $\widetilde{\mathfrak G}^r_d$ restricting to the already
defined vector bundle $\H_{a, b}$ over $\sigma^{-1}(\cM_g^0\cup
\Delta_0^0\cup \Delta_1^0)$, such that $\H_{0,
b}=\rm{Sym}$$^b(\G_{0, 1})$ for all $b\geq 1$ and which also has the
property that the exact sequences (\ref{hi}) remain exact over
$\widetilde{\mathfrak G}^r_d$.
\item For each $0\leq a\leq r, b\geq 1$, there exists a torsion free sheaf
 $\G_{a, b}$ over $\widetilde{\mathfrak G}^r_d$ that restricts to the
vector bundle $\G_{a, b}$ over $\sigma^{-1}(\cM_g^0\cup \Delta_0^0
\cup \Delta_1^0)$, which for $a=0$ agrees with the vector bundle
$\G_{0, b}$ defined above, and which has the property that the
vector bundle morphisms $\phi_{a, b}$ defined over
$\sigma^{-1}(\cM_g^0\cup \Delta_0^0\cup \Delta_1^0)$ extend to
morphisms $\phi_{a, b}:\H_{a, b}\rightarrow \G_{a, b}$ over
$\widetilde{\mathfrak G}^r_d$.
\end{enumerate}
\end{proposition}
\begin{proof} We start with an arbitrary point $t=(C\cup_y D, l_C, l_D)\in
\sigma^*(C^j)$ where we assume first that $[g/2]\leq j\leq rs-s$. We
set $L_C:=l_C(-(rs-j)y)\in \mbox{Pic}^{r+j}(C)$ and
$L_D:=l_D(-(r+j)y)\in \mbox{Pic}^{rs-j}(D)$. If $L=(L_C, L_D)\in
\mbox{Pic}^{rs+r}(C\cup_y D)$, the essential observation is that
because $[C]\in \cM_j$ and $[D, y]\in \cM_{g-j, 1}$ are
Brill-Noether general, we always have that $rs+r-1\leq
a^{l_C}_0(y)+a^{l_D}_r(y)\leq rs+r$, hence $h^0(L_C)\geq r$,
$h^0(L_D)\leq 1$ and $h^0(C\cup_y D, L)=r+1$ (see Proposition
\ref{limitlinj}). If $p:\widetilde{\cM}_{g, 1}\times
_{\widetilde{\cM}_g}\widetilde{\mathfrak G}^r_d \rightarrow
\widetilde{\mathfrak G}^r_d$ is the universal curve over
$\widetilde{\mathfrak G}^r_d$, we denote by $\mathcal{P}$ a
Poincar\'e bundle of relative degree $d=rs+r$ enjoying the following
properties: 

(1) For each $[g/2]\leq j\leq g-2$, $\mathcal{P}_{|
p^{-1}(\sigma^{-1}(\Delta_j^0))}$ parameterizes line bundles of
bidegree $(r+j, rs-j)$ on curves of type $C\cup_y D$ where $g(C)=j$
and  $g(D)=g-j$. \ 

(2) If $\tau_j:\sigma^{-1}(\Delta_j^0)\rightarrow
\widetilde{\cM}_{g, 1}\times_{\widetilde{\cM}_g}
\widetilde{\mathfrak G}^r_d$ denotes the section which assigns the
single node corresponding to every curve from $\sigma^{-1}(\Delta_j^0)$,
then $\tau_j^*(\mathcal{P})=\OO_{\sigma^{-1}(\Delta_j^0)}$. 

Note that since the divisors $\sigma^{-1}(\Delta_i^0)$ and
$\sigma^{-1}(\Delta_j^0)$ are disjoint for $[g/2]\leq i<j\leq g-2$,
the construction can be carried out over a fixed divisor
$\Delta_j^0$ at a time. Since $h^0(p^{-1}(t), \mathcal{P}_{|
p^{-1}(t)})=h^0(C\cup_y D, L)=r+1$ for each $t\in
\widetilde{\mathfrak G}^r_d$,  by Grauert's Theorem, $\G_{0,
1}:=p_*(\mathcal{P})$ is a locally free sheaf which satisfies our
first requirement. For $b\geq 2$ we define $\G_{0,
b}=p_*(\mathcal{P}^{\otimes b})$. Based on degree considerations we
have that $H^1(L_C^{\otimes b})=H^1(L_C^{\otimes b}\otimes
\OO_C(-y))=0$. Using Proposition \ref{vanish2} we see that
$H^1(L_D^{\otimes b})=0$, hence $h^0(C\cup D, L^{\otimes
b})=h^0(L_C^{\otimes b})+h^0(L_D^{\otimes b})-1=bd+1-g$. Grauert's
Theorem implies that $\G_{0, b}$ is locally free for all $b$.

In the remaining case when $rs-s+1\leq j\leq g-2$, that is, $2\leq
g(D)\leq s-1$, we define $L_C:=l_C(-(rs+s+r-j)y) \in
\mbox{Pic}^{j-s}(C)$ and $L_D:=l_D(-(j-s)y)\in
\mbox{Pic}^{rs+s-j+r}(D)$.  Then $h^0(L_C)\leq 1, h^0(L_D)\geq r$
and $h^0(C\cup_y D, L)=r+1$. Proposition \ref{vanish2} gives again
that $h^1(L_C^{\otimes b})=0$ for all $b\geq 2$,  hence $h^0(C\cup
D, L^{\otimes b})=bd+1-g$. This time , we denote by $\mathcal{P}$
the Poincar\'e bundle parametrizing line bundles of bidegree $(j-s,
rs+s+r-j)$ on curves of type $C\cup_y D$ and then $\G_{0,
b}:=p_*(\mathcal{P}^{\otimes b})$ is locally free in this case too
because of Grauert's Theorem.

To define $\G_{a, b}$ for $a\geq 1$, we introduce the sheaf
$\cM:=\mbox{Ker}\{p^*(p_*(\mathcal{P}))\rightarrow \mathcal{P}\}$
and then we set $\G_{a, b}:=p_*(\wedge^a \cM\otimes
\mathcal{P}^{\otimes b})$. The morphism $\phi_{0, b}$ is simply the
natural map $\mbox{Sym}^b p_*(\mathcal{P})\rightarrow
p_*(\mathcal{P}^{\otimes b})$, and to define these maps for $a\geq
1$ we use that the vector bundles $\H_{a, b}$ fit into exact
sequences of type (\ref{hi}) and then proceed inductively.
\end{proof}
\begin{remark} Just like in the case of the vector bundles $\cA$ and $\cB$ defined initially over $\mathfrak{G}^r_d$ (cf. Section 2), the sheaves $\H_{a, b}, \G_{a, b}$ depend on the choice of a Poincar\'e bundle, whereas $Hom\ _{\OO_{\overline{\mathfrak G}^r_d}} (\H_{a, b}, \G_{a, b})$ and $\phi_{a, b}\in H^0\bigl(\overline{\mathfrak G}^r_d, Hom\ _{\OO_{\overline{\mathfrak G}^r_d}} (\H_{a, b}, \G_{a, b})\bigr)$ are independent of such a choice. Moreover, since the projection $p:\widetilde{\cM}_{g, 1}\times_{\widetilde{\cM}_g} \widetilde{\mathfrak G}^r_d\rightarrow \widetilde{\mathfrak G}^r_{d}$ has a canonical section over each divisor $\sigma^{-1}(\Delta_j^0)$ where $[g/2]\leq j\leq g-1$, it is possible to choose the Poincar\'e bundle $\mathcal{P}_{| p^{-1}(\sigma^{-1}(\Delta_j^0))}$ in an unambiguous way (which is precisely what we did in the proof of Proposition \ref{ext0}) and then $\H_{a, b| \sigma^{-1}(\Delta_{j}^0)}$ and $\G_{a, b | \sigma^{-1}(\Delta_j^0)}$ are unambiguously defined as vector bundles over $\sigma^{-1}(\Delta_j^0)$. This is of course  a minor point which plays no role in the calculation of  $\sigma_*(c_1(\G_{i, 2}-\H_{i, 2}))\in A^1(\widetilde{\cM}_g)$.
\end{remark}

We now determine the class of the curves $X$ and $Y$ defined in
Proposition \ref{limitlin}:

\begin{proposition}\label{xy}
Let $C$ be a Brill-Noether general curve of genus $g-1$ and $q\in C$
a general point. We denote by $\pi_2: C\times W^r_d(C)\rightarrow
W^r_d(C)$ the projection and set
$c_i:=(\pi_2)^*\bigl(c_i(\mathcal{E}^{\vee})\bigr)$.

\noindent(1) The class of the curve $X=\{(y,L)\in C\times
W^r_d(C):h^0(C,L(-2y))\geq r\}$ is given by
$$[X]= c_r+c_{r-1}(2\gamma+(2d+2g-4)\eta)-6c_{r-2}\  \eta \theta.
$$

\noindent (2) The class of the curve $Y=\{(y,L)\in C\times
W^r_d(C):h^0(C,L(-y-q))\geq r\}$ is given by
$$[Y]= c_r+c_{r-1}(\gamma +(d-1)\eta)-2c_{r-2}\ \eta \theta.$$
\end{proposition}
\begin{proof}
We realize both $X$ and $Y$ as degeneracy loci over $C\times
W^r_d(C)$ and compute their classes using the Thom-Porteous formula.
For each $(y, L)\in C\times W^r_d(C)$ we have a natural map $$
H^0(C, L_{|2y})^{\vee} \rightarrow H^0(C, L)^{\vee}$$ which
globalizes to a vector bundle map $\zeta: J_1(\mathcal{L})^{\vee}
\rightarrow (\pi_2)^*(\mathcal{E}^{\vee})$.  Clearly
$X=Z_{1}(\zeta)$, hence
$$[X]=\bigl[\frac{c_t(\pi_2^*(\mathcal{E}^{\vee}))}{c_t(J_1(\mathcal{L})^{\vee})}\bigr]_r\ .$$
\noindent From the exact sequence defining the jet bundle of $\mathcal{L}$
$$0\longrightarrow \pi_1^*(K_C)\otimes \mathcal{L}\rightarrow
J_1(\mathcal{L})\rightarrow \mathcal{L}\rightarrow 0$$ we obtain
that
$c_t(J_1(\mathcal{L}^{\vee}))^{-1}=1+2\gamma+2d\eta+(2g-4)\eta-6\eta
\theta$, which quickly leads to the desired expression for $[X]$.
The calculation of $[Y]$ is entirely similar and we skip it.
\end{proof}

We also need the following intersection theoretic result:
\begin{lemma}\label{fj}
For each $j\geq 2$ we have the following formulas:
\begin{enumerate}
\item
$c_1(\G_{0, j \ |X})=-j^2\theta-(2g-4)\eta-j(d\eta+\gamma)$.
\item
$c_1(\G_{0, j \ |Y})=-j^2\theta+\eta$.
\end{enumerate}
\end{lemma}
\begin{proof} We observe that for all $j\geq 2$, $H^1(L^{\otimes j})=0$, hence
 $(\pi_2)_*(\mathcal{L}^{\otimes j})$ is a vector bundle over
$\mbox{Pic}^d(C)$.
 Riemann-Roch applied to the
map $\pi_2:C\times \mbox{Pic}^d(C)\rightarrow \mbox{Pic}^d(C)$
yields $c_1\bigl((\pi_2)_*(\mathcal{L}^{\otimes
j})\bigr)=-j^2\theta$. If we denote by $u, v:C\times C\times
\mbox{Pic}^d(C)\rightarrow C\times \mbox{Pic}^d(C)$ the two
projections and by $\Delta\subset C\times C\times \mbox{Pic}^d(C)$
the diagonal, we have the following exact sequences
$$0\longrightarrow u_*(v^*(\mathcal{L}^{\otimes j})\otimes
\mathcal{I}_{\Delta}^2)_{| X} \longrightarrow \G_{0, j | X} \longrightarrow
\mathcal{L}^{\otimes j}_{|X}\longrightarrow 0$$ and
$$0\longrightarrow u_*(v^*(\mathcal{L}^{\otimes j})\otimes
\mathcal{I}_{\Delta}^2)\longrightarrow
(\pi_2)_*(\mathcal{L}^{\otimes j})\longrightarrow
J_1(\mathcal{L}^{\otimes j})\longrightarrow 0$$ (and an entirely
similar situation for $\G_{0, j |Y}$) from which both claims follow
easily.
\end{proof}

Now we are in a position to prove Theorem \ref{div}:
\newline
\emph{Proof of Theorem \ref{div}.}
 Since
$\mbox{codim}(\mm_g-\widetilde{\cM}_g, \mm_g)\geq 2$, it makes no
difference whether we compute the class $\sigma_*(\G_{i, 2}-\H_{i,
2})$ on $\widetilde{\cM}_g$ or on $\mm_g$ and we can write
\begin{equation}\label{classd}
\sigma_*(\G_{i, 2}-\H_{i, 2})= A \lambda-B_0\ \delta_0-B_1\
\delta_1-\cdots -B_{[g/2]}\ \delta_{[g/2]}\in \mbox{Pic}(\mm_g),
\end{equation}
where $\lambda, \delta_0, \ldots, \delta_{[g/2]}$ are the generators
of $\mbox{Pic}(\mm_g)$. We start with the following:

\noindent {\bf Claim:} \emph{One has the relation
$A-12B_0+B_1=0$.}

We pick a general curve $[C, q]\in \cM_{g-1, 1}$ and at the fixed
point $q$ we attach to $C$ a Lefschetz pencil of plane cubics. If we
denote by $R\subset \mm_g$ the resulting curve, then
 $R\cdot \lambda=1, \ R\cdot \delta_0=12, \ R\cdot
\delta_1=-1$ and $R\cdot \delta_j=0$ for $j\geq 2$. The relation
$A-12B_0+B_1=0$ follows once we show that $\sigma^*(R)\cdot
c_1(\G_{i, 2}-\H_{i, 2})=0$. To achieve this we check that $\G_{0, b
|\sigma^*(R)}$ is trivial and then use (\ref{gi}) and (\ref{hi}). We
take $[C\cup_q E]\in \mm_g$ to be an arbitrary curve from $R$, where
$E$ is an elliptic curve. The pointed curve $[C, q]$ being
Brill-Noether general, limit $\mathfrak g^{r}_{d}$'s on $C\cup_q E$
are in $1:1$ correspondence with linear series $L\in W^{r}_{d}(C)$
having a cusp at $q$ (This is a statement independent of the
$j$-invariant of $E$, in particular, it also holds for the $12$
rational nodal curves in the pencil). Furthermore, the fibre of
$\G_{0, b | \sigma^*(\Delta_1^0)}$ over each point from
$\sigma^*(R)$ consists of the global sections of the genus $g-1$
aspect of the limit $\mathfrak g^{r}_{d}$ and the claim now follows.

Now we determine explicitly the coefficients $A, B_0$ and $B_1$. We
fix a general curve $[C, q]\in \cM_{g-1, 1}$ and construct the test
curves $C^1\subset \Delta_1$ and $C^0\subset \Delta_0$. Using the
notation from Proposition \ref{limitlin}, we write that
$\sigma^*(C^0)\cdot c_1(\G_{i, 2}-\H_{i, 2}) =c_1(\G_{i, 2
|Y})-c_1(\H_{i, 2 |Y})$   and
 $\sigma^*(C^1) \cdot c_1(\G_{i, 2}-\H_{i, 2})=c_1(\G_{i, 2 |X})-c_1(\H_{i, 2 |X})$
 (the other component $T$ of $\sigma^*(C^1)$ does not appear
 because $\G_{0, b |T}$ is trivial for all $b\geq 1$).
On the other hand $$C^0 \cdot \sigma_*(c_1(\G_{i, 2}-\H_{i,
2}))=(2g-2)B_0-B_1 \mbox{ and } C^1\cdot \sigma_*(c_1(\G_{i,
2}-\H_{i, 2}))=(2g-4)B_1,$$ while we already know that
$A-12B_0+B_1=0.$ Next we use the relations
$$c_1(\G_{i, 2})= \sum_{l=0}^i (-1)^l c_1(\wedge^{i-l}\G_{0,
1}\otimes \G_{0, l+2})=\sum_{l=0}^i (-1)^l{ r+1 \choose i-l}
c_1(\G_{0, l+2})+$$ $$+\sum_{l=0}^i (-1)^l
\bigl((l+2)(rs+r)+1-rs-s\bigr){r \choose i-l-1}c_1(\G_{0, 1}), \
\mbox{ } \mbox{ and }$$
$$c_1(\H_{i, 2})=\sum_{l=0}^i (-1)^l c_1(\wedge^{i-l} \G_{0,
1}\otimes \mbox{Sym}^{l+2} \G_{0, 1})=$$ $$=\sum_{l=0}^i
(-1)^l\Bigl( {r \choose i-l-1}{r+l+2 \choose l+2}+{r+1\choose
i-l}{r+l+2\choose r+1}\Bigr)c_1(\G_{0, 1})=$$
$$={2s+is+i\choose i} (s+1)(i+2)c_1(\G_{0, 1}),$$ which when restricted
to $X$ and $Y$, enable us (also using Lemma \ref{fj}), to obtain
explicit expressions for $c_1(\G_{i, 2}-\H_{i, 2})_{|X}$ and
$c_1(\G_{i, 2}-\H_{i, 2})_{| Y}$ in terms of the classes $\eta,
\theta, \gamma$ and $c_1=\pi_2^*(c_1(\mathcal{E}^{\vee}))$.
Intersecting these classes with $[X]$ and $[Y]$ and using Lemma
\ref{vandermonde}, we finally get a linear system of $3$ equations
in $A, B_0$ and $B_1$ which leads to the stated formulas for the
first three coefficients. \hfill $\Box$

Next we prove that when $i=0$, we can get a formula for the slope of
$\zz_{s(2s+1), 0}$: precisely we show that if we write
$\sigma_*(c_1(\G_{0, 2}-\H_{0, 2}))=A\lambda-B_0\delta_0-\cdots
-B_{[g/2]}\delta_{[g/2]},$ then $B_j\geq B_0$ for all $j\geq 1$. In
particular, $s\bigl(\sigma_*(c_1(\G_{0, 2}-\H_{0, 2}))\bigr)=A/B_0$
which has already been computed in Theorem \ref{div}. We note that
the proof uses in an essential way the divisor class calculation
from Theorem \ref{pointedbn}.

\noindent \emph{Proof of Theorem \ref{khosla}.} Using the convention
$B_{g-j}=B_j$ for $g/2\leq j\leq g-1$, we show that $B_j\geq B_0$
only when $s(2s+1)/2\leq j\leq s(2s-1)$. The case $2s^2\leq j\leq
s(2s+1)-1$ is dealt with in a similar fashion. To compute $B_j$ we
intersect the class $\sigma_*(c_1(\G_{0, 2}-\H_{0, 2}))$ with $C^j$.
Then we use that $[Y_{g-j, \beta}]\cdot c_1(\G_{0, 2}-\H_{0,
2})=[U_{j, \gamma}]\cdot c_1(\G_{0, 2}-\H_{0, 2})=0$, for all
$\beta\in \mathcal{P}_2, \gamma\in \mathcal{P}_3$, to obtain that
$$(2j-2)B_j=\sigma^*(C^j)\cdot c_1(\G_{0, 2}-\H_{0,
2})=\sum_{(\alpha_0, \ldots, \alpha_r)\in \mathcal{P}_1} N_{g-j,
\alpha} \bigl([X_{j, \alpha}]\cdot c_1(\G_{0, 2}-\H_{0, 2})\bigr).$$
We fix a Schubert index $(\alpha_0\leq \ldots \leq \alpha_r)\in
\mathcal{P}_1$ and denote by $\pi_1:X_{j, \alpha}\rightarrow C$ and
$\pi_2:X_{j, \alpha}\rightarrow \mbox{Pic}^{r+j}(C)$ the two
projection maps. As before, $\mathcal{L}$ is the Poincar\'e bundle
on $C\times \mbox{Pic}^{r+j}(C)$. There is an isomorphism of bundles
$\G_{0, 1| X_{j,
\alpha}}=\pi_2^*\bigl((\pi_2)_*(\mathcal{L})\bigr)_{| X_{j,
\alpha}}$ obtained by globalizing the projection isomorphism at the
level of spaces of sections $H^0(C\cup_y D, L)\cong H^0(C, L_C)$
valid for each point $(y, L_C)\in X_{j, \alpha}$. (We recall that
$L=(L_C, L_D)\in \mbox{Pic}^{2s+j}(C)\times
\mbox{Pic}^{2s^2-j}(D)$). For $b\geq 2$, we have a surjective
morphism of vector bundles $\G_{0, b | X_{j,
\alpha}}\twoheadrightarrow
\pi_2^*\bigl((\pi_2)_*(\mathcal{L}^{\otimes b})\bigr)_{| X_{j,
\alpha}}$ whose kernel is a trivial bundle along $X_{j, \alpha}$.
Thus one has that $c_1(\G_{0, b |X_{j, \alpha}})=-b^2 \theta_{|
X_{j, \alpha}}$ and $c_1(\H_{0, 2 |X_{j, \alpha}})=c_1(\mbox{Sym}^2
\G_{0, 1| X_{j, \alpha}})=-(2s+2)\theta_{| X_{j, \alpha}}$,
therefore
\begin{equation}\label{estim}
(2j-2)B_j=(2s-2)\sum_{(\alpha_0, \ldots, \alpha_r)\in \mathcal{P}_1}
N_{g-j, \alpha} \ ([X_{j, \alpha}]\cdot \theta).
\end{equation}
 The class of the curve $X_{j, \alpha}$ can be computed using the generalized
 Giambelli formula (cf. \cite{FuPr}, pg. 15-17) as follows: If
 $J_{\alpha_r+r-1}(\mathcal{L})\twoheadrightarrow \cdots \twoheadrightarrow
 J_{\alpha_i+i-1}(\mathcal{L})\twoheadrightarrow \cdots \twoheadrightarrow J_{\alpha_0-1}(\mathcal{L})$
 is the flag of jet bundles  corresponding to the ramification
 sequence $(\alpha_0, \ldots, \alpha_r)$, then $$X_{j, \alpha}=\{(y, L)\in C\times \mbox{Pic}^{r+j}(C):
 \mbox{rk}\{\pi_2^*\bigl((\pi_2)_*(\mathcal{L})\bigr)(y, L)\rightarrow J_{\alpha_i+i-1}(\mathcal{L})(y, L)\}\leq i
 \mbox{ for all }i\}
 $$
 and then $[X_{j, \alpha}]$ is given by the determinant of the $(\alpha_r\times
 \alpha_r)$-matrix having entries
 $$a_{ik}=c_{r+1-l+k-i}\Bigl(\frac{\pi_2^*\bigl((\pi_2)_*(\L)\bigr)}{J_{\alpha_{l}+l-1}(\L)}\Bigr),
 \mbox{ for all } \alpha_{l-1}\leq i\leq \alpha_l, 0\leq l\leq r \mbox{ and } 1\leq j\leq
 \alpha_r.$$
 Since
 $$c_i\Bigl(\frac{\pi_2^*\bigl((\pi_2)_*(\L))}{J_{a-1}(\L)}\Bigr)=\frac{\theta^i}{i!}+\frac{\theta^{i-1}}{(i-1)!}\bigl(a\gamma+
 a(j+r)\eta+a(a-1)(j-1)\eta \bigr)-\frac{\theta^{i-1}}{(i-2)!}a(a+1)\eta,$$
 clearly $[X_{j, \alpha}]$
 is a linear combination of $\theta^j, \theta^{j-1}\eta$ and
 $\theta^{j-1}\gamma$ in $H^{2j}(C\times \mbox{Pic}^{r+j}(C))$.
The intersection number $[X_{j, \alpha}]\cdot \theta$ can be
interpreted as the number of line bundles $L_C\in
\mbox{Pic}^{r+j}(C)$ satisfying the condition $\alpha^{L_C}_i(y)\geq
\alpha_i$ for $i=0, \ldots, r$ at an \emph{unspecified} point $y\in
C$, and which, moreover, are also ramified at a \emph{fixed} point
$q\in C$, that is, $a_r^{L_C}(q)\geq r+1$.

Using this interpretation, the quantity $\sum_{\alpha \in
\mathcal{P}_1} N_{g-j, \alpha} ([X_{g-j, \alpha}]\cdot \theta)$ can
be expressed as the intersection number $\widetilde{C}_j\cdot
\overline{\mathfrak{Lin}}^r_d(1)$ over the moduli space $\mm_{g,
1}$.  Here $\mathfrak{Lin}^r_d(1)$ is the divisor on $\cM_{g, 1}$
consisting of pointed curves $[C, q]$ such that there exists $L\in
W^r_d(C)$ with $h^0(C, L\otimes \OO_C(-(r+1)q)\geq 1$, while
$\widetilde{C}_j=\{C\cup_y D, q\}_{y\in C} \subset
\Delta_{j:1}\subset \mm_{g, 1}$ is the test curve obtained by
varying the point of attachment $y$ on the genus $j$ component,
while the marked point $q\in C$ remains fixed. The class of
$\overline{\mathfrak{Lin}}^r_d(1)$ is computed in the course of the
proof of Theorem \ref{pointedbn} and one has
$$\overline{\mathfrak{Lin}}^r_d(1)\equiv
\mu\Bigl((g+3)\lambda-\frac{g+1}{6}\delta_{irr}-\sum_{j=1}^{g-1}\delta_{j:1}\Bigr)+
\nu\Bigl(-\lambda+\psi-\sum_{j=1}^{g-1} {g-j+1\choose
2}\delta_{j:1}\Bigr),$$ where
$$\nu=\frac{r(r+2)}{(rs+s-1)(rs+s+1)} \mbox{ and
}\mu=\frac{r(r+1)(r+2)(s-1)(s+1)(rs+s+4)}{2(s+r+1)(rs+s-2)(rs+s-1)(rs+s+1)}.$$
Since $\widetilde{C}_j\cdot \psi=1, \widetilde{C}_j\cdot
\delta_{g-j: 1}=1$ (the only point of intersection corresponds to
$y=q \in C$), $\widetilde{C}_j \cdot \delta_{j:1}=-(2j-1)$, while
$\widetilde{C}_j\cdot \lambda=\widetilde{C}_j \cdot \delta_{i:1}=0$
for $i\neq j, g-j$, we can compute that
$$\frac{B_j}{c_r}=\frac{s-1}{(j-1)c_r}\sum_{\alpha \in
\mathcal{P}_1}N_{g-j, \alpha}([X_{j, \alpha}]\cdot
\theta)=\frac{s-1}{(j-1)c_r}\ \widetilde{C}_j\cdot
\overline{\mathfrak{Lin}}^r_d(1)=$$
$$=\frac{4(s-1)j\bigl(2js^3+js^2-2js-2j+4s^3+4s^2-3s\bigr)(2s^2+s-j)}{(2s^2+s-2)(3s+1)(2s-1)(j-1)}\geq$$
$$\geq \frac{B_0}{c_r}=
\frac{s(8s^6-8s^5-2s^4+s^2+11s+2)}{3(2s^2+s-2)(3s+1)(2s-1)}\ .$$
This finishes the proof and shows that $s(\sigma_*(\G_{0, 2}-\H_{0,
2}))=A/B_0$. \hfill $\Box$

As we have already pointed out, Theorem  \ref{div} produces only
virtual divisors on $\mm_g$ of slope less than $6+12/(g+1)$. To get
actual divisors one has to show  that the vector bundle map
$\phi:\H_{i, 2}\rightarrow \G_{i, 2}$ is generically non-degenerate.
We carry this out in the case $i=0$ and we produce for the first
time an infinite sequence of genuine counterexamples to the Slope
Conjecture.

\noindent \emph{Proof of \ Theorem  \ref{maxrank}.} From
Brill-Noether theory one knows that there exists a unique component
of $\widetilde{\mathfrak G}^r_d$ which maps onto
$\widetilde{\cM}_g$. Moreover, if $(C, L)\in \mathfrak G^r_d$ is
such that $L\in W^r_d(C)-W^{r+1}_d(C)$ corresponds to an embedding
$C\subset \PP^r$, then a sufficient condition for the smoothness of
$\mathfrak{G}^r_d$ at $[C, L]$ is that $H^1(N_{C/\PP^r})=0$, and
then, the differential $(d \sigma)_{[C, L]}$ is surjective if and
only if the \emph{Petri map} $\mu_0(C): H^0(L)\otimes H^0(K_C\otimes
L^{\vee})\rightarrow H^0(K_C)$ is injective (see e.g. \cite{AC2}).
In our situation, it is then enough to produce a Brill-Noether-Petri
general smooth curve $C\subset \PP^{2s}$ having degree $2s(s+1)$ and
genus $s(2s+1)$ such that $C$ does not sit on any quadrics, that is
$H^0(\mathcal{I}_{C/\PP^{2s}}(2))=H^1(\mathcal{I}_{C/\PP^{2s}}(2))=0$.
We carry this out inductively: For each $0\leq a\leq s$, we
construct a smooth non-degenerate curve $C_a \subset \PP^{s+a}$ with
$\mathrm{deg}(C_a)={s+a+1 \choose 2}+a$ and $g(C_a)={s+a+1\choose
2}+a-s$, $h^1(C_a, \OO_{C_a}(1))=a$ (or equivalently, $h^0(C_a,
\OO_{C_a}(1))=s+a+1$),  and such that (1) $C_a$ satisfies the Petri
Theorem (in particular one has that $H^1(C_a, N_{C_a
/\PP^{s+a}})=0$), and (2) the multiplication map $\mu_2:
\mbox{Sym}^2 H^0(C_a, \OO_{C_a}(1))\rightarrow H^0(C_a,
\OO_{C_a}(2))$ is surjective (or equivalently, an isomorphism).

To construct $C_0\subset \PP^s$ we consider the \emph{White surface}
$S=\mbox{Bl}_{\{p_1, \ldots, p_{\delta}\}}(\PP^2) \subset \PP^s$
obtained by blowing-up $\PP^2$ at general points $p_1, \ldots,
p_{\delta}\in \PP^2$ where $\delta={s+1\choose 2}$,  and embedding
it via the linear system $|s h-\sum_{i=1}^{\delta} E_{p_i}|$. Here
$h$ is the class of a line on $\PP^2$. It is known that $S\subset
\PP^s$ is a projectively Cohen-Macaulay surface and its ideal is
generated by the $(3\times 3)$-minors of a certain $(3\times
s)$-matrix of linear forms (see e.g. \cite{GG} even though these
surfaces have been studied in the classical literature by T.G. Room
in \cite{R}). The Betti diagram of $S\subset \PP^s$ is the same as
that of the ideal of $(3\times 3)$-minors of a $(3\times s)$-matrix
of indeterminates. In particular, we have that
$H^i(\mathcal{I}_{S/\PP^s}(2))=0$ for $i=0, 1$. On $S$ we consider a
generic smooth curve $C \equiv (s+1)h-\sum_{i=1}^{\delta} E_{p_i}$.
We find that the embedded curve $C\subset S\subset \PP^s$ has
$\mbox{deg}(C)={s+1\choose 2}$ and $g(C)={s \choose 2}$. From the
exact sequence $$0\longrightarrow
\mathcal{I}_{S/\PP^s}(1)\longrightarrow \mathcal{I}_{C/\PP^s}(1)
\longrightarrow  \mathcal{I}_{C/S}(1)\longrightarrow 0,$$ using also
that $H^1(\mathcal{I}_{S/\PP^s}(1))=0$ and that
$H^1(\mathcal{I}_{C/S}(1))=0$ (e.g. by Riemann-Roch), we find that
$H^1(\mathcal{I}_{C/\PP^s}(1))=0$ and $H^1(\OO_{C}(1))=0$, hence
$h^0(\OO_{C}(1))=s+1$. Furthermore, since
$H^0(\mathcal{I}_{S/\PP^s}(2))=H^1(\mathcal{I}_{S/\PP^s}(2))=0$, we
obtain that $H^1(\mathcal{I}_{C/\PP^s}(2))=0$. Finally, since
$H^1(\OO_{C}(1))=0$, it follows trivially that $H^1(N_{C/\PP^s})=0$
and $\mu_0(C)$ is injective, being a map with source the trivial
vector space. Even though $[C]\in \cM_{g(C)}$ itself is not a Petri
general curve, the map $\mathcal{H}_C\rightarrow \cM_{g(C)}$ from
the Hilbert scheme  $\mathcal{H}_C$ of curves $C'\subset \PP^s$ with
$\mbox{deg}(C')=\mbox{deg}(C)$ and $g(C')=g(C)$, is smooth and
dominant around the point $[C]\in \mathcal{H}_C$, hence a generic
deformation $[C_0\hookrightarrow \PP^s]\in \mathcal{H}_C$ of
$[C\hookrightarrow \PP^s]$ will be Petri general and still satisfy
the condition $H^1(\mathcal{I}_{C_0/\PP^s}(2))=0$.

Assume now that we have constructed a Petri general curve
$C_a\subset \PP^{s+a}$ with all the desired properties. We pick
general points $p_1, \ldots, p_{s+a+2}\in C_a$ with the property
that if $\Delta:=p_1+\cdots +p_{s+a+2}\in \mbox{Sym}^{s+a+2}C_{a}$,
then the variety
$$T:=\{(M, V)\in W^{s+a+1}_{d(C_a)+s+a+2}(C_a):
\mbox{dim}\bigl(V\cap H^0(C_a, M\otimes
\OO_{C_a}(-\Delta))\bigr)\geq s+a+1\}$$ of linear series having an
$(s+a+2)$-fold point along $\Delta$, has the expected dimension
$\rho(g(C_a), s+a+1, d(C_a)+s+a+2)-(s+a+1)^2$. We identify the
projective space $\PP^{s+a}$ containing $C_a$ with a hyperplane
$H\subset \PP^{s+a+1}$ and choose a linearly normal elliptic curve
$E\subset \PP^{s+a+1}$ such that $E\cap H=\{p_1, \ldots,
p_{s+a+2}\}$. The fact that such an $E$ exists is an easy
consequence of the vanishing $H^1(N_{E/\PP^{s+a+1}}(-1))=0$ for each
elliptic curve $E$ embedded by a complete linear series; the
vanishing itself is a consequence of the fact that
$N_{E/\PP^{s+a+1}}$ is a poly-stable vector bundle (cf. \cite{GL2},
Theorem 4.1), having the property that
$\mu(N_{E/\PP^{s+a+1}}(-1))>1$). We now set $X:=C_a\cup_{\{p_1,
\ldots, p_{s+a+2}\}}E\hookrightarrow \PP^{s+a+1}$ and then
$\mbox{deg}(X)=p_a(X)+s$. From the exact sequence $$0\longrightarrow
\OO_E(-p_1-\cdots -p_{s+a+2})\longrightarrow \OO_X\longrightarrow
\OO_{C_a}\longrightarrow 0,$$ we can write that $h^0(\OO_X(1))\leq
h^0(\OO_{C_a}(1))+h^0(\OO_E)=s+a+2$, hence $h^0(\OO_X(1))=s+a+2$ and
$h^1(\OO_X(1))=a+1$. One can also write the exact sequence
$$0\longrightarrow \mathcal{I}_{E/\PP^{s+a+1}}(1)\longrightarrow
\mathcal{I}_{X/\PP^{s+a+1}}(2)\longrightarrow
\mathcal{I}_{C_a/H}(2)\longrightarrow 0,$$ from which we obtain that
$H^1(\mathcal{I}_{X/\PP^{s+a+1}}(2))=0$, hence by a dimension count
we also get that $H^0(\mathcal{I}_{X/\PP^{s+a+1}}(2))=0$, that is,
$X$ and a general deformation of $X$ inside $\PP^{s+a+1}$  lie on no
quadrics.

We now show that $X\hookrightarrow \PP^{s+a+1}$ can be deformed to
an embedding of a smooth curve $C_{a+1}$ in $\PP^{s+a+1}$ such that
$H^1(N_{C_{a+1}/\PP^{s+a+1}})=0$. We choose an $(s+a+2)$-dimensional
subspace $H^0(\OO_{C_a}(1))\subset V\subset H^0(\OO_{C_a}(1)\otimes
\OO_{C_a}(\Delta))$ which gives a map $f: C_a\rightarrow
\PP^{s+a+1}$ such that $f(p_1)=\cdots =f(p_{s+a+2})=p$. If we denote
by $\widetilde{\PP}^{s+a+1}$ the blow-up of $\PP^{s+a+1}$ at $p$, by
choosing $V$ suitably we may assume that $f$ lifts to an embedding
$\tilde{f}:C_a\hookrightarrow \widetilde{\PP}^{s+a+1}$ which
projected from $p$ gives rise to the original embedding
$C_a\hookrightarrow H$. We consider another copy of $\PP^{s+a+1}$
which we denote by $\PP^{s+a+1}_1$ and we denote by $Z$ the scheme
obtained by gluing $\PP^{s+a+1}_1$ and $\widetilde{\PP}^{s+a+1}$
along $H$, where we identify the exceptional divisor of
$\widetilde{\PP}^{s+a+1}$ with $H\subset \PP^{s+a+1}$ via the
projection from $p$. There is a natural map $h: Z\rightarrow
\PP^{s+a+1}$ which on $\PP^{s+a+1}_1$ is the identity while on
$\widetilde{\PP}^{s+a+1}$ is the projection from $p$. Via this map,
the inclusion $X\hookrightarrow \PP^{s+a+1}$ lifts to an embedding
$X \hookrightarrow Z$. Note that $Z$ is a degeneration of
$\PP^{s+a+1}$ something which can be seen by blowing-up the
codimension $2$ subscheme $H\times \{0\}$ of $\PP^{s+a+1}\times
\PP^1$. If we denote by $\mathcal{X}$ the total space of the blow-up
and by $\epsilon: \mathcal{X}\rightarrow \PP^1$ the projection onto
the second component, then for $t\neq 0$ we have that
$\epsilon^{-1}(t)=\PP^{s+a+1}$, whereas $\epsilon^{-1}(0)=\PP\cup
\E$, where $\PP$ is the strict transform of $\PP^{s+a+1}\times
\{0\}$ which is  isomorphic to $\PP^{s+a+1}$, while
$\E=\PP(\OO_H\oplus \OO_H(1))$ is the exceptional divisor, which is
isomorphic to $\PP^{s+a+1}$ blown-up at a point. In the special
fibre, $\PP$ and $\E$ are joined along a divisor which is $H$ inside
$\PP$.

Next we write down the standard exact sequences of normal bundles
$$0\longrightarrow N_{E/\PP^{s+a+1}}\otimes \OO_{E}(-\Delta)\longrightarrow
N_{X/Z}\longrightarrow
N_{C_a/\widetilde{\PP}^{s+a+1}}\longrightarrow 0 $$ (the right hand
side map is restriction to the component $C_a$ of $X$), \   and
$$0\longrightarrow \OO_{C_a}(1)\otimes \OO_{C_a}(2 \Delta) \longrightarrow N_{C_a/\widetilde{\PP}^{s+a+1}}
\longrightarrow N_{C_a/H} \longrightarrow 0,$$ from which it easily
follows that $H^1(N_{X/Z})=0$  (Use the hypothesis
$H^1(N_{C_a/\PP^{s+a}})=0$ and that
$H^1\bigl(N_{E/\PP^{s+a+1}}\otimes \OO_{E}(-\Delta)\bigr)=0$ because
$N_{E/\PP^{s+a+1}}$ is semi-stable). Thus the space of deformations
of $X$ in $Z$ is unobstructed of dimension $h^0(N_{X/Z})$. On the
other hand, by general theory the space $T^1_{(X, Z)}$ of
infinitesimal deformations of the pair $(X, Z)$ has dimension at
least $\chi(N_{X/Z})+1=h^0(N_{X/Z})+1$. This shows that there exists
a deformation of $(X, Z)$ in which $Z$ deforms non-trivially. But
$\mbox{dim} (T_Z^1)=1$, that is, the only possible deformation of
$Z$ is the smoothing to $\PP^{s+a+1}$ previously described, and in
this deformation the map $X\hookrightarrow Z$ will deform to an
embedding $C_{a+1}\hookrightarrow \PP^{s+a+1}$ of a smooth curve,
which proves our claim. We are left with showing that the dimension
estimate
\begin{equation}\label{dimest}
\mbox{dim}\bigl(W^{s+a+1}_{d(C_{a+1})}(C_{a+1})\bigr)=\rho\bigl(g(C_{a+1}),
s+a+1, d(C_{a+1})\bigr)
\end{equation} holds. Assuming that (\ref{dimest}) has been proved,
since the  condition $H^1(N_{C_{a+1}/\PP^{s+a+1}})=0$ guarantees the
local smoothness of the scheme $\mathfrak G^{s+a+1}_{d(C_{a+1})}$,
it follows that the morphism $\mathfrak
G^{s+a+1}_{d(C_{a+1})}\rightarrow \cM_{g(C_{a+1})}$ is dominant in a
neighbourhood of the point $[C_{a+1}\hookrightarrow \PP^{s+a+1}]$.
Therefore the curve $C_{a+1}\subset \PP^{s+a+1}$ can be chosen to be
Petri general as well, which enables us to continue the induction.

We return to proving (\ref{dimest}) and denote by $\mathcal{U}$ the
versal deformation space of $[X]\in \mm_{g(C_{a+1})}$ and by $\phi:
\mathcal{C}\rightarrow \mathcal{U}$ the universal family such that
$\phi^{-1}(0)=X$, where $0\in \cU$. Then in a way similar to
\cite{EH1}, Theorem 3.3, one can construct a quasi-projective
variety $\sigma:\tilde{\mathfrak G}_{d(C_{a+1})}^{s+a+1}\rightarrow
\mathcal{U}$ of limit linear series such that for points $u\in
\mathcal{U}$ with $C_u=\phi^{-1}(u)$ smooth, we have that
$\sigma^{-1}(u)=G^{s+a+1}_{d(C_{a+1})}(C_u)$, whereas
$\sigma^{-1}(0)$ consists of the following data: an underlying line
bundle $\mathcal{L}$ on $X$ together with linear series $\{L_{a},
V_{a}\in G(s+a+2, H^0(X, L_{a}))\}$ and $\{L_E, V_E\in G(s+a+2,
H^0(X, L_E))\}$ such that the following conditions are satisfied
(see also \cite{Est}, Theorem 1):
\begin{enumerate}
\item The line bundles
$L_{a}$ and $L_E$ on $X$ are suitable twists of $\mathcal{L}$ by
multiples of the divisor $\Delta$: precisely there exists an integer
$l$ such that $L_{a |C_a}=L_{E| C_a}\otimes \OO_{C_a}(l\cdot
\Delta)$ and $L_{a |E}=L_{E | E}\otimes \OO_E(-l\cdot \Delta)$.
Moreover
$\mathrm{deg}(L_{a|C_a})+\mathrm{deg}(L_{E|E})=\mathrm{deg}(C_{a+1})+l(s+a+2)$.
\item The restriction maps $V_{a}\rightarrow H^0(C_{a}, L_{a
|C_{a}})$ and $V_E\rightarrow H^0(E, L_{E | E})$ are both injective.
\item
The restriction maps $V_{a}\rightarrow H^0(E, L_{a |E})$ and
$V_E\rightarrow H^0(C_{a}, L_{E |C_{a}})$ are both non-zero.
\item If $l$ is the integer defined above and $(a_0\leq \ldots \leq a_{s+a+1})$ denotes the vanishing sequence of
$(L_{a| C_{a}}, V_{a+1})$ with respect to the divisor $\Delta \in
\mbox{Sym}^{s+a+2} C_{a}$ while $(b_0\leq \cdots \leq b_{s+a+1})$
denotes the vanishing sequence of $(L_{E | E}, V_E)$ with respect to
$\Delta \in \mbox{Sym}^{s+a+2} E$, then we have the inequalities
$a_i+b_{s+a+1-i}\geq l$ for all indices $0\leq i\leq s+a+1$  (see
also \cite{Est}, Proposition 6).
\end{enumerate}
By construction we have the dimension estimate
$\mbox{dim}(\mathcal{G})\geq \mbox{dim}(\cU)+\rho\bigl(g(C_{a+1}),
s+a+1, d(C_{a+1})\bigr)$, (see also \cite{EH1}), thus in order to
prove (\ref{dimest}) it suffices to show that
$$\mbox{dim}(\sigma^{-1}(0))=\rho\bigl(g(C_{a+1}), s+a+1,
d(C_{a+1})\bigr)=\rho(g(C_a), s+a, d(C_a))-a$$ (here by dimension we
mean the smallest dimension of an irreducible component). It is now
easy to describe the fibre $\sigma^{-1}(0)$ in a neighbourhood of
the point corresponding to a smoothing of the embedding
$X\hookrightarrow \PP^{s+a+1}$: If $\{L_a, V_a\}, \{L_E, V_E\}$ is a
limit linear series on $X$, then the aspect corresponding to $E$ is
just a very ample line bundle on $X$ giving the embedding into
$\PP^{s+a+1}$, that is, $L_{E|C_a}\in W^{s+a}_{d(C_a)}(C_a)$ and
$L_{E|E}=\OO_E(\Delta)$, whereas the aspect corresponding to $C_a$
is described by $L_{a|E}=\OO_E$ and $L_{a
|C_a}=\OO_{C_a}(\Delta)\otimes L_{E | C_{a}}$ (and in particular
$l=1$). The only possibility for the vanishing sequences of the $E$
and $C_a$ aspects is that $(a_0, \ldots,  a_{s+a+1})=(0, 1, \ldots
,1)$ and $(b_0, \ldots, b_{s+a+1})=(0, \ldots, 0, 1)$. This shows
that locally, $\sigma^{-1}(0)$ is isomorphic to the variety of line
bundles $\L\in \mbox{Pic}^{d(C_{a+1})}(X)$ such that
$\L_{|E}=\OO_E$, $h^0(X,\L)\geq s+a+2$ and $h^0(C_a,
\L_{|C_a}(-\Delta))\geq s+a+1$ (Loosely speaking this is the
subscheme consisting of those $L_a\in W^{s+a}_{d(C_a)}(C_a)$ for
which there exists a section $\tau \in \PP\bigl(H^0(L_a\otimes
\OO(\Delta))/H^0(L_a)\bigr)$ which glues to the unique section of
the trivial bundle $\OO_E$ at the points of attachment $p_1, \ldots,
p_{s+a+2}$). Thus locally $\sigma^{-1}(0)$ is a $(\mathbb
C^*)^{s+a+1}$-bundle over the subvariety $T$ of
$W_{d(C_{a+1})}^{s+a+1}$ having an $(s+a+2)$-fold point along the
divisor $\Delta$, and by our inductive hypothesis we know that
$\mbox{dim}(T)=\mbox{dim}\bigl(
W_{d(C_{a+1})}^{s+a+1}(C_a)\bigr)-(s+a+1)^2$. It follows that
$$\mbox{dim}(\sigma^{-1}(0))=\mbox{dim}(T)+s+a+1=\rho\bigl(g(C_{a+1}), s+a+1, d(C_{a+1})\bigr),$$
and this finishes the proof. \hfill
$\Box$

\begin{remark} It is natural to wonder whether (\ref{dimest})
could not be proved more easily by showing directly that the Petri
map $\mu_0(X): H^0(\OO_{X}(1))\otimes H^0(\omega_X(-1))\rightarrow
H^0(\omega_X)$ is injective. Indeed Theorem 1.3 from \cite{CR2}
seems to imply this to be the case based on the inductive hypothesis
that $\mu_0(C_a)$ is injective whereas $\mu_0(E)$ is injective for
trivial reasons. That claim is incorrect: from the exact sequence
$0\longrightarrow \omega_E\longrightarrow \omega_X \longrightarrow
\omega_{C_a}(\Delta)\longrightarrow 0$, we find the isomorphism
$H^0(\omega_X(-1))=H^0(\omega_{C_a}(-1)\otimes \OO_{C_a}(\Delta))$
and then a simple analysis shows that $H^0(\OO_E)\otimes
H^0(\omega_{C_a}(-1))\subset H^0(\OO_X(1))\otimes H^0(\omega_X(-1))$
is an $a$-dimensional subspace lying entirely inside
$\mbox{Ker}(\mu_0(X))$.
\end{remark}

\section{The class of the Gieseker-Petri divisors}

In this section we prove Theorem \ref{gp}. We use the same strategy
as in the previous section and we intersect $\GisP$ with the test
curves $C^0, C^1$ and $C^j$ for $[g/2]\leq j \leq g-2$. Recall that
we have constructed a rank $r+1$ vector bundle $\G_{0, 1}$ over the
variety $\widetilde{\mathfrak{G}}^r_d$ (cf. Proposition \ref{ext0}).
As usual, we denote by $\mathbb E$ the Hodge bundle over $\mm_g$.

\begin{proposition}
There exists a rank $s$ vector bundle $\cN$ over
$\widetilde{\mathfrak{G}}^r_d$ together with a morphism $\G_{0,
1}\otimes \cN\rightarrow \sigma^*(\mathbb E\otimes
\OO_{\widetilde{\cM}_g}(\delta_1))$ of vector bundles of the same
rank over $\widetilde{\mathfrak{G}}^r_d$ such that the fibres of
$\cN$ admit the following description:
\begin{itemize}
\item If $(C, L)\in \mathfrak{G}^r_d$, then $\cN(C, L)=H^0(C,
K_C\otimes L^{\vee})$.
\item If $t=(C\cup_y E, L_C, l_E)\in \sigma^{-1}(\Delta_1^0)$, where
$L\in W^r_d(C)$ is such that $h^0(L\otimes \OO_C(-2y))=r$, then
$\cN(t)=H^0(C, K_C\otimes L_C^{\vee}\otimes \OO_C(2y))$.
\item If $t=(C/y\sim q, L)\in \sigma^{-1}(\Delta_0^0)$, where $y,
q\in C$ and $L\in W^r_d(C)$ is a linear series such that
$h^0(L\otimes \OO_C(-y-q))=h^0(L)-1$, then $\cN(t)=H^0(C, K_C\otimes
L^{\vee}\otimes \OO_C(y+q))$.
\item If $t=(C\cup_y D, l_C, l_D)\in \sigma^{-1}(\Delta_j^0)$ where
$[g/2]\leq j\leq g-2, g(C)=j, g(D)=g-j$,  then $\cN(t)=H^0(C\cup_y
D, \omega_{C\cup D} \otimes L^{\vee})$, where $L=(l_C(-(rs-j)y),
l_D(-(j+r)y)) \in \rm{Pic}$$^{j+r}(C)\times \rm{Pic}$$^{rs-j}(D)$.
\end{itemize}
\end{proposition}

Note that over $\mathfrak{G}^r_d$ the morphism $\G_{0, 1}\otimes
\cN\rightarrow \sigma^*(\mathbb E\otimes
\OO_{\widetilde{\cM}_g}(\delta_1))$ is simply the Petri
multiplication map. We start the proof of Theorem \ref{gp} by
expanding $[\GisP]$ in $\mbox{Pic}(\mm_g)$:
$$\GisP\equiv  a\lambda-b_0\delta_0-\cdots
-b_{[g/2]}\delta_{[g/2]}.$$ We show that the coefficients $a, b_0$
and $b_1$ as well as $s(\GisP)$ can be read off from the vector
bundle map $\G_{0, 1}\otimes \cN\rightarrow \sigma^*(\mathbb
E\otimes \OO_{\widetilde{\cM}_g}(\delta_1))$. Even though this
bundle map is degenerate along the boundary components contained in
$\sigma^*(\Delta_j^0)$ with $j\geq 2$, we can show that it is
generically non-degenerate along $\sigma^*(\Delta_0^0)$ and
$\sigma^*(\Delta_1^0)$ which ultimately suffices to compute
$s(\GisP)$.

\begin{proposition}\label{gispell}
One has the relation $a-12b_0+b_1=0$. Moreover, one has the identity
$$\GisP\equiv \sigma_*\Bigl(c_1(\sigma^*(\mathbb E\otimes \OO_{\widetilde{\cM}_g}(\delta_1))-c_1(\G_{0,
1}\otimes \cN)\Bigr)+\sum_{j=2}^{[g/2]} d_j \delta_j, $$ where
$d_j\geq 0$.
\end{proposition}
\begin{proof} It is enough to show that if $[C, y]\in \cM_{g-1, 1}$ is a
general pointed curve, then for every $L\in W^r_d(C)$ satisfying
$h^0(L\otimes \OO_C(-2y))=r$, the multiplication  map
$$\mu_{0}(L, y):
H^0(L)\otimes H^0(K_C\otimes L^{\vee}\otimes \OO_C(2y))\rightarrow
H^0(K_C\otimes \OO_C(2y))$$ is an isomorphism. This shows that the
morphism $\G_{0, 1}\otimes \cN\rightarrow \sigma^*(\mathbb E\otimes
\OO_{\widetilde{\cM}_g}(\delta_1))$ is non-degenerate along each
component of the divisor $\sigma^{-1}(\Delta_1^0)$ and the
conclusion follows. To show that $\mu_0(L, y)$ is an isomorphism, we
use a variation of the degeneration considered by Eisenbud and
Harris to prove the Gieseker-Petri Theorem (cf. \cite{EH4}).
Precisely, we consider a $1$-dimensional family
$\pi:\mathcal{C}\rightarrow B$ of generically smooth pointed curves
of genus $g-1$  with a section $\tau: B\rightarrow \mathcal{C}$,
degenerating to a curve of compact type $C_0$ consisting of a string
of rational components and $g$ elliptic components $E_1, \ldots,
E_g$ such that the stable model of $C_0$ is $E_1\cup_{p_1}
E_2\cup_{p_2} E_3\cup \cdots \cup_{p_{g-1}} E_{g-1}$. We assume
moreover that the marked point specializes to a point $p_0\in E_1$
and we choose our degeneration general enough such that
$p_i-p_{i-1}\in \mbox{Pic}^0(E_i)$ is not a torsion point for all
$1\leq i\leq g-1$. By contradiction, we assume that for a general
$t\in B$ there exists $L_t\in W^r_d(\pi^{-1}(t))$ with
$h^0\bigl(\pi^{-1}(t), L_t\otimes \OO(-2\tau(t))\bigr)=r$, such that
$\mu_0(L_t, \tau(t))$ has non-trivial kernel. For $1\leq i \leq g$
we denote by $L^i\in \mbox{Pic}^d(C_0)$ the limit line bundle of the
$L_t$'s having the property that $\mbox{deg}(L^i_{| E_j})=0$ for
$i\neq j$, hence $\mbox{deg}(L^i_{| E_i})=d$. Similarly, we define
 $M^i \in \mbox{Pic}^{2g-2-d}(C_0)$ to be the limit when $t\rightarrow
0$ of $K_{\pi^{-1}(t)}\otimes L_t^{\vee}\otimes
\OO_{\pi^{-1}(t)}(2\tau(t))$ uniquely characterized by the property
$\mbox{deg}(M^i_{| E_j})=0$ for $i\neq j$ and $\mbox{deg}(M^i_{|
E_i})=2g-2-d$. We denote by $\{\bigl(L^i_{| E_i}, V_i\in G(r+1,
H^0(E_i, L_{|E_i}))\bigr)\}$ and by $\{\bigl(M^i_{| E_i}, W_i \in
G(r+1, H^0(E_i, M_{| E_i}))\bigr)\}$ the limit linear series on
$C_0$ corresponding to $L_t$ and $K_{\pi^{-1}(t)}\otimes L_t^{\vee}$
respectively as $t\rightarrow 0$. Reasoning along the lines of
\cite{EH4} or \cite{F2}, Proposition 5.2, for each $1\leq i\leq g$
we find non-trivial elements
$$\rho_i\in \mbox{Ker}\{V_i\otimes W_i\rightarrow H^0(E_i,
L^i\otimes M^i_{| E_i})\}$$ satisfying
$\mbox{ord}_{p_i}(\rho_{i+1})\geq \mbox{ord}_{p_{i-1}}(\rho_{i})+2$
for $1\leq i\leq g-1$. Since both $V_1$ and $W_1$ have a cusp at
$p_0\in E_1$, it follows that $\mbox{ord}_{p_1}(\rho_1)\geq 2$,
hence $\mbox{ord}_{p_{g-1}}(\rho_g)\geq
2g-2=\mbox{deg}(L^g_{|E_g})+\mbox{deg}(M^g_{| E_g})$, which is a
contradiction because $\rho_{g-1}\in H^0(E_{g}, L^g_{| E_g})\otimes
H^0(E_g, M^g_{| E_g})$ being an element in the kernel of the
multiplication map must be a tensor of rank at least $4$.
\end{proof}

\begin{proposition}\label{gispd1}
If $c_r$ is the constant defined in Lemma \ref{vandermonde}, then
the $\delta_1$ coefficient in the expression of $[\GisP]$ is given
by:
$$ b_1=c_r \frac{r(s-1)}{(s+r+1)(rs+s-2)}(3rs^2+2s^2+r^2s^2+7s+6rs+r^2s+2r+2).$$
\end{proposition}
\begin{proof} We fix a general curve $C$ of genus $g-1$ and
consider the associated test curve $C^1\subset \Delta_1$. We view
the curve $X\subset C\times W^r_d(C)$ defined in Proposition
\ref{xy}, as sitting inside $\widetilde{\mathfrak G}^r_d$. Then the
projection $\pi_1:X \rightarrow C$ is the restriction of
$\sigma:\widetilde{\mathfrak G} ^r_d \rightarrow \widetilde{\cM}_g$
once we identify $C$ with $C_1$ (Note that the degree of $\pi_1$ is
precisely $c_r$). One can write the relation $(2g-4)B_1=C^1 \cdot
\GisP=c_1(\sigma^*(\mathbb E\otimes
\OO_{\widetilde{\cM}_g}(\delta_1))_{\ |X})-c_1(\mathcal{G}_{0, 1 \
|X}\otimes \mathcal{N}_{\ |X})$ and we are going to compute each
term in this expression.

The restriction  $\mathbb E\otimes
\OO_{\widetilde{\cM}_g}(\delta_1)_{|C^1}$ is identified with the
vector bundle $(p_2)_*\bigl(p_1^*(K_C)\otimes \OO(2\Delta)\bigr)$,
where $p_1, p_2:C\times C\rightarrow C$ are the two projections and
$\Delta\subset C\times C$ is the diagonal. Using
Grothendieck-Riemann-Roch for the map $p_2$, we find that
$$c_1(\mathbb E\otimes \OO_{\widetilde{\cM}_g}(\delta_1)_{|  C^1})=c_1\bigl((p_2)_{!}(p_1^*(K_C)\otimes
\OO(2\Delta))\bigr)=-2g+4,$$ hence $c_1(\sigma^*(\mathbb E)\otimes
\OO_{\widetilde{\cM}_g}(\delta_1))_{| X}=-(2g-4)c_r$ (remember that
we have set $c_i=c_1(\mathcal{E}^{\vee})$).

The fibre $\mathcal{N}_{| X}(y, L)$ is identified with
$H^0(K_C\otimes L^{\vee}(2y))=H^1(L\otimes \OO(-2y))^{\vee}$.
Keeping in mind that we have introduced the vector bundle map
$\zeta$ in Proposition \ref{xy}, we have an exact sequence over $X$
$$0\longrightarrow \mbox{Ker}(\zeta)^{\vee}\longrightarrow
\mathcal{N}_{ |X}^{\vee}\longrightarrow \pi_2^*\bigl(R^1\pi_{2
*}(\mathcal{L}_{| C\times W^r_d(C)})\bigr)\longrightarrow 0,$$
globalizing the cohomology exact sequence for each $(y, L)\in
C\times W^r_d(C)$
$$ \cdots \longrightarrow
H^0(L)\stackrel{\zeta^{\vee}}\longrightarrow H^0(L_{|
2y})\longrightarrow H^1(L(-2y))\longrightarrow H^1(L)\longrightarrow
0.$$ Hence $c_1(\mathcal{N}_{
|X}^{\vee})=\theta-c_1(\mathcal{E}^{\vee})+c_1(\mbox{Ker}(\zeta)^{\vee})$.
Using Proposition \ref{xy}  we can write that
$$C^1\cdot \GisP=-(2g-4)c_r \eta -(r+1-s)c_1\cdot
[X]+(r+1)c_1(\mbox{Ker}(\zeta)^{\vee})=$$
$$
-(2g-4)c_r \eta -(r+1-s)\bigl((2d+2g-4)c_1
c_{r-1}\eta-6c_1c_{r-2}\theta \eta\bigr)+$$
$$+(r+1)\bigl((2d+2g-4)c_{r-1}\theta \eta-6c_{r-2}\theta^2
\eta\bigr)+(r+1)c_1(\mbox{Ker}(\zeta)^{\vee}).$$ To compute the
Chern number $c_1(\mbox{Ker}(\zeta)^{\vee})$ we use once more
\cite{HT} and we find the following relation in $H^{top}(C\times
W^r_d(C))$:
$$c_1(\mbox{Ker}(\zeta)^{\vee})=c_{r+1}\Bigl(\frac{\pi_2^*(\mathcal{E}^{\vee})}{J_1(\mathcal{L}^{\vee})
}\Bigr)=(2d+2g-4)c_r\eta-6\eta \theta c_{r-1}.$$ Combining the last
two relations and then applying Lemma \ref{vandermonde} we obtain
the formula for $b_1$.
\end{proof}

\begin{proposition}\label{gispd0}

The $\delta_0$ coefficient in the expression of $[\GisP]$ is given
by:
$$b_0=c_r\frac{r(r+1)(r+2)(s-1)s(s+1)(rs+s+4)}{6(r+s+1)(rs+s-2)(rs+s-1)}.$$
\end{proposition}
\begin{proof} We pick a general curve $C$ of genus $g-1$ and consider the test
curve $C^0\subset \Delta_0$. Similarly to the proof of Proposition
\ref{gispd1} we view the projection $\pi_1:Y \rightarrow C$ as the
restriction of $\sigma:\widetilde{\mathfrak G}^r_d \rightarrow
\widetilde{\cM}_g$ over $C^0$. Then one has the relation
$$(2g-2)b_0-b_1=C^0\cdot \GisP=c_1(\sigma^*(\mathbb E\otimes \OO_{\widetilde{\cM}_g}(\delta_1)
)_{| Y})-c_1(\mathcal{G}_{0, 1|Y}\otimes \mathcal{N}_{|Y}).$$ The
Hodge bundle $\mathbb E\otimes \OO_{\widetilde{M}_g}(\delta_1)_{|
C^0}$ is identified with $(p_2)_*\bigl(p_1^*(K_C)\otimes
\OO(\Delta+\Gamma_q)\bigr)$, where $\Gamma_q=\{q\}\times C$, and it
is easy to compute that $c_1(\sigma^*(\mathbb E\otimes
\OO_{\widetilde{\cM}_g}(\delta_1))_{|Y})=c_r$. If we denote by
$\upsilon$ the vector bundle morphism over $Y$ which globalizes the
maps $H^0(L_{| y+q})^{\vee}\rightarrow H^0(L)^{\vee}$ for each $(y,
L)\in Y$, we obtain an exact sequence of vector bundles over $Y$
$$0 \longrightarrow \mbox{Ker}(\upsilon)^{\vee} \longrightarrow \mathcal{N}_{| Y}^{\vee}\longrightarrow
\pi_{2}^*\bigl(R^1\pi_{2 *}(\mathcal{L}_{| C\times
W^r_d(C)})\bigr)\longrightarrow 0,$$ from which we can compute
$c_1(\mathcal{N}_{| Y}^{\vee})$ if we use \cite{HT} which in this
case reads $$c_1(\mbox{Ker}(\upsilon
^{\vee}))=c_{r+1}\Bigl(\frac{\pi_2^*(\mathcal{E}^{\vee})}{\mathcal{F}}\Bigr),$$
where $\mathcal{F}$ is the vector bundle on $C\times W^r_d(C)$ with
fibre $\mathcal{F}(y, L)=H^0(L_{| y+q})^{\vee}$. Finally, we write
$$C^0\cdot \GisP=c_r\eta+((r+1)(\theta-c_1+c_1))\cdot [Y]+(r+1)(c_r(d-1)-2c_{r-1}\theta)\eta,$$
which eventually leads to the stated formula.
\end{proof}
To finish the proof of Theorem \ref{gp} it suffices to show that for
$[g/2]\leq j\leq g-2$, the coefficient of $\delta_j$ in the
expression of $\sigma_*\bigl(c_1(\sigma^*(\mathbb E)\otimes
\OO_{\widetilde{\cM}_g}(\delta_1))-c_1(\G_{0, 1}\otimes \cN)\bigr)$
always exceeds the coefficient of $\delta_0$, which equals $b_0$ and
was computed in Proposition \ref{gispd0}. This is a calculation
along the lines of the proof of Theorem \ref{khosla}. To keep the
length of this paper under control, we skip the details.

\section{Five ways of constructing Koszul divisors for pointed curves}

In this section we construct Koszul divisors on moduli spaces
of pointed curves. As an application we improve Logan's results on
which $\mm_{g,n}$'s are of general type.

We start by recalling a few things about divisor classes on
$\mm_{g,n}$.
 For $0\leq i\leq g$ and a set of indices $S\subset \{1,\ldots,n\}$, the boundary divisor
 $\Delta_{i:S}$  corresponds to the closure of the locus of nodal curves
 $C_1\cup C_2$, with $C_1$ smooth of genus $i$, $C_2$ smooth
 of genus $g-i$, and such that the marked points sitting on
 $C_1$ are precisely those labeled by $S$.
We also introduce the divisor  $\Delta_{irr}$ consisting of
irreducible pointed curves with one node. We denote by
$\delta_{i:S}\in \mbox{Pic}(\overline{\mathcal{M}}_{g,n})$ the class
of $\Delta_{i:S}$ and by $\delta_{irr}$ that of $\Delta_{irr}$. For
each $1\leq i \leq n$
 we define the tautological class
$\psi_i:=c_1(\mathbb L_i)$, where $\mathbb L
_i$ is the line bundle
on $\mm_{g,n}$ with fibre $\mathbb L_i([C, x_1, \ldots,
x_n])=T_{x_i}^{\vee}(C)$ over each point $[C, x_1, \ldots, x_n]\in \mm_{g, n}$.
 It is well known that when $g\geq 3$, the Hodge class $\lambda$, the boundaries $\delta_{irr}$ and $\delta_{i:S}$,
 and the tautological classes $\psi_i$ for $1\leq i\leq n$, freely generate
 $\mbox{Pic}(\overline{\mathcal{M}}_{g,n})$.

\subsection{Divisors defined in terms of the Minimal Resolution
Conjecture}\hfill

We fix integers $g, r\geq 1$ and $0\leq i\leq g$ and set
$n:=(2r+1)(g-1)-2i$. We define a divisor on $\cM_{g, n}$ consisting
of smooth pointed curves $(C, x_1, \ldots, x_n)$ such that the
points $x_1, \ldots, x_n$ fail the \emph{Minimal Resolution
Conjecture} for the canonical curve
$C\stackrel{|K_C|}\hookrightarrow \PP^{g-1}$ (see \cite{FMP} for
background on MRC). Precisely we define the locus
$$\mathfrak{Mrc}_{g,i}^r:=\{[C,x_1, \ldots, x_n]\in \cM_{g,n}: h^1\bigl(C, \wedge^{i}
M_{K_C}\otimes K_C^{\otimes (r+1)}\otimes
\OO_C(-x_1-\cdots-x_n)\bigr)\geq 1\}.$$

If we denote by $\Gamma:=x_1+\cdots+x_n\in C_n$, by Serre duality,
the condition appearing in the definition of $\mathfrak{Mrc}_{g,
i}^r$ is equivalent to
$$h^0\bigl(C, \wedge^{i} M_{K_C}^{\vee}\otimes
\OO_C(\Gamma)\otimes K_C^{\otimes (-r)}\bigr)\geq 1
\Longleftrightarrow \OO_C(\Gamma)\otimes K_C^{\otimes (-r)}\in
\Theta_{\wedge^{i} M_{K_C}^{\vee}},$$ where we recall that for a
stable vector bundle $E$ on $C$ having  slope $\nu(E)=\nu\in \mathbb
Z$, its \emph{theta divisor} is the determinantal locus
$$\Theta_E:=\{ \eta \in \mbox{Pic}^{g-\mu-1}(C): h^0(C, E\otimes
\eta) \geq 1 \}.$$ The main result from \cite{FMP} gives an
identification $\Theta_{\wedge^i M_{K_C}^{\vee}}=C_{g-i-1}-C_i$,
where the right hand side is one of the difference varieties
associated to $C$. Thus one has an alternative description of points
in $\mathfrak{Mrc}_{g,i}^r$: a point $(C, x_1, \ldots, x_n)\in
\mathfrak{Mrc}_{g, i}^r$ if and only if there exists $D\in
C_i$ such that $h^0\bigl(C, \OO_C(\Gamma+D)\otimes K_C^{\otimes
(-r)}\bigr)\geq 1$.

First we equip $\mathfrak{Mrc}_{g, i}^r$ with a determinantal scheme
structure. We consider the following cartesian diagram of stacks
\[
\begin{CD}
{\mathcal X}@>{q}>> \cM_{g,n}\\
@VV{f}V@VV{\pi}V\\
{\mathcal C}_g^{}@>{p}>>\cM_{g}\\
\end{CD}
\]
in which all the morphisms are smooth and $p$ (hence also $q$) is
proper. We denote by $\omega_{p} \in \mbox{Pic}(\mathcal{C}_g)$ the
relative dualizing sheaf of the universal curve
$p:\mathcal{C}_g\rightarrow \mathcal{M}_g$ and by $\mathbb
E:=p_*(\omega_p)$ the Hodge bundle. We define the vector bundle
$\mathcal{M}$ over $\mathcal{C}_g$ having rank $g-1$ as the kernel
of the evaluation map $p^*\mathbb E\longrightarrow\omega_p$. Thus
for every $[C]\in\cM_g$, we have $\cM \vert_{p^{-1}([C])}\simeq
M_{K_C}$. For each $1\leq j\leq n$ we have a section
$q_j\,:\,\cM_{g,n}\longrightarrow {\mathcal X}$ of $q$ given by
$q_j([C, x_1, \ldots, x_n])=([C, x_1, \ldots, x_n], x_j)\in
\mathcal{X}$ and we set $E_j:={\rm Im}(q_j)$, hence $E_j$ is a
relative divisor over $\cM_{g,n}$.

For integers $0\leq a\leq  i, b\geq r+2$ and $(a, b)=(0, r+1)$ we
define the vector bundle
$$\cA_{a,b}:=q_*\bigl(f^*(\wedge^a \cM\otimes \omega_p^{\otimes b})\otimes
\OO_{\mathcal{X}}(-\sum_{j=1}^n E_j)\bigr),$$ hence $\cA_{a, b}([C,
x_1, \ldots, x_n])=H^0(\wedge^a M_{K_C}\otimes K_C^{\otimes
b}\otimes \OO_C(-\Gamma))$. To prove that $\cA_{a, b}$ is locally
free over $\cM_{g, n}$, we use the fact that $M_{K_C}$ is a
semi-stable vector bundle over $C$ and that $\mu\bigl(\wedge^a
M_{K_C}\otimes K_C^{\otimes b}(-\Gamma)\bigr)>2g-1$, hence
$$R^1q_*\bigl(f^*(\wedge^a \cM \otimes \omega_p^{\otimes b})\otimes
\OO_{\mathcal{X}}(-\sum_{j=1}^n E_j)\bigr)=0$$ whenever $b\geq r+2$.
To reach the same conclusion in the case of the sheaf $\cA_{0,
r+1}$, we use that $H^1(K_C^{\otimes
(r+1)}(-\Gamma))^{\vee}=H^0(\OO_C(\Gamma)\otimes K_C^{\otimes
(-r)})=0$, if $\Gamma \in C_n$ lies outside a subset of codimension
$\geq 3$.  Furthermore there is a vector bundle map $$\phi:
\pi^*(\wedge^i \mathbb E)\otimes \cA_{0, r+1}\rightarrow \cA_{i-1,
r+2}$$ which over each point $[C, x_1, \ldots, x_n]\in \cM_{g, n}$
corresponds to the natural map
$$\phi(C, \Gamma): \wedge^i H^0(K_C)\otimes H^0\bigl(K_C^{\otimes (r+1)}\otimes
\OO_C(-\Gamma)\bigr)\rightarrow H^0\bigl(\wedge^{i-1} M_{K_C}\otimes
K_C^{\otimes (r+2)}\otimes \OO_C(-\Gamma)\bigr).$$ Note that
$\mbox{rank}(\cA_{i-1, r+2})=\mbox{rank}\bigl(\pi^*(\wedge^i \mathbb
E)\otimes \cA_{0, r+1}\bigr)=2i{g \choose i}$ and a simple argument
using the exact sequence $0\longrightarrow \wedge^i
M_{K_C}\longrightarrow \wedge^i H^0(K_C)\otimes \OO_C\longrightarrow
\wedge^{i-1} M_{K_C}\otimes K_C\longrightarrow 0$ shows that
$\mathfrak{Mrc}_{g, i}^r$ is the degeneracy locus of the map $\phi$.

\begin{proposition}\label{hyper}
The vector bundle morphism $\phi: \pi^*(\wedge^i \mathbb E)\otimes
\cA_{0, r+1}\rightarrow \cA_{i-1, r+2}$ is generically
non-degenerate. It follows that $\mathfrak{Mrc}_{g, i}^r$ is a
divisor on $\cM_{g, n}$.
\end{proposition}

\begin{proof} We show that $\phi$ is generically non-degenerate over
the pull-back $\pi^*(\mathcal{H}_g)$ of the hyperelliptic locus. We
fix a hyperelliptic curve $C$ of genus $g$ and we denote by $L\in
W^1_2(C)$ its hyperelliptic involution. By writing down the Euler
sequence on $\PP^1$ one shows that $M_{K_C}=(L^{\vee})^{\oplus
(g-1)}$, hence the condition $H^1(\wedge^i M_{K_C}\otimes
K_C^{\otimes (r+1)}\otimes \OO_C(-\Gamma))=0$ is equivalent to
$H^0\bigl(\OO_C(\Gamma )\otimes L^{\otimes  (i-r(g-1))}\bigr)=0$.
This however is obvious because when $\Gamma \in C_n$ is a general
divisor of degree $n$ then $\OO_C(\Gamma)\otimes L^{\otimes
(i-r(g-1))}$ is a general line bundle of degree $g-1$, therefore it
has no global sections.
\end{proof}

The main result here is the computation of the class of
$\overline{\mathfrak{Mrc}}_{g, i}^r$:

\begin{theorem}\label{mrc}
When $n=(2r+1)(g-1)-2i$, the locus $\mathfrak{Mrc}_{g,i}^r$ is a
divisor on $\cM_{g,n}$ and the class of its compactification in
$\mm_{g, n}$ is given by the following formula:
$$\overline{\mathfrak{Mrc}}_{g,i}^r\equiv \frac{1}{g-1}{g-1\choose
i}\Bigl(a\lambda+c\sum_{j=1}^n \psi_j-b_{irr}\delta_{irr}-\sum_{j,
s\geq 0, } b_{j:s}\sum_{|S|=s} \delta_{j:S}\Bigr),$$ where
$$c=rg+g-i-r-1, \ b_{irr}=-\frac{1}{g-2}\Bigl({r+1\choose 2}(g-1)(g-2)+i(i+1+2r-rg-g)\Bigr),$$
$$a=-\frac{1}{g-2}\Bigl((g-1)(g-2)(6r^2+6r+1)+i(24r+10i+10-10g-12rg)\Bigr),$$
$$
\ b_{0:s}={s+1\choose 2}(g-1)+s(rg-r)-si, \mbox{ and } b_{j:s}\geq
b_{0: s} \mbox{ for } j\geq 1.$$
\end{theorem}
\begin{remark} For $i=0, n=(2r+1)(g-1)$, the
divisor $\overline{\mathfrak{Mrc}}_{g, 0}^r$ specializes to the
locus of points $[C, x_1, \ldots, x_{(2r+1)(g-1)}]\in \cM_{g,
(2r+1)(g-1)}$ such that $\sum_{j=1}^{(2r+1)(g-1)} x_j\in
|K_C^{\otimes r}|$ and Theorem \ref{mrc} gives that:
$$\overline{\mathfrak{Mrc}}_{g, 0}^r\equiv
-(6r^2+6r+1)\lambda+(r+1)\Bigl(\sum_{j=1}^{(2r+1)(g-1)}
\psi_j\Bigr)+ {r+1\choose
2}\delta_{irr}-(2r+3)\sum_{|S|=2}\delta_{0:S}-\cdots.$$ By letting
all the marked points coalesce and using the standard formulas for
pushing forward products of tautological classes (cf. e.g.
\cite{FMP} or \cite{Log}, Theorem 2.8), Theorem \ref{mrc} offers a
quick way of computing the class of the closure of the locus
$\mathcal{W}^{r+1}$ of $(r+1)$-Weierstrass points in $\mm_{g, 1}$
which is the main result of \cite{CF}.
\end{remark}

First we determine the class of the locus $\mathfrak{Mrc}_{g, i}^r$
over the interior $\cM_{g, n}$. In order to do this, we first recall
a few well-known intersection theory relations (see e.g. \cite{HM}):

\begin{lemma}\label{rules}
If $q:\mathcal{X}\rightarrow \cM_{g, n}$ is the morphism defined
earlier, one has the following identities:
\item {\rm (i)} $q_*(f^*c_1(\omega_p)^2)=12\lambda$.
\item {\rm (ii)} $q_*(q^*\lambda\cdot f^*c_1(\omega_p))=(2g-2)\lambda$.
\item {\rm (iii)} $q_*q^*(\lambda^2)=0$.
\item {\rm (iv)} $q_*(c_1(E_j)\cdot q^*\lambda)=\lambda$.
\item {\rm (v)} $q_*(c_1(E_j)\cdot f^*c_1(\omega_p))=\psi_j$.
\item {\rm (vi)} $q_*q^*c_2(\pi^*(\mathbb E))=0$.
\item {\rm (vii)} $q_*(c_1(E_j)^2)=-\psi_j$.
\end{lemma}

\begin{proposition}\label{interior}
If $a$ and $b$ are the numbers defined in the statement of Theorem
\ref{mrc}, we have the following relation in $\mbox{Pic}(\cM_{g,
n})$:
$$\mathfrak{Mrc}_{g, i}^r\equiv \frac{1}{g-1}{g-1 \choose i}\bigl(a\lambda+b\sum_{j=1}^n
\psi_j\bigr).$$
\end{proposition}
\begin{proof} Since $\phi$ is generically non-degenerate (cf.
Proposition \ref{hyper}), we have the identity $\mathfrak{Mrc}_{g,
i}^r\equiv c_1(\cA_{i-1,r+2})-c_1(\pi^*(\wedge^i \mathbb E)\otimes
\cA_{0, r+1}).$ To compute these Chern classes we use
Grothendieck-Riemann-Roch applied to the proper map $q$. For
simplicity we set $\mathfrak D:=\sum_{j=1}^n E_j$ and
$\mathcal{F}:=f^*\bigl(\wedge^{i-1} \mathcal{M}\otimes
\omega_p^{\otimes (r+2)})\otimes \OO_{\mathcal{X}}(-\mathfrak
D)\bigr)$. Then we have that
$$c_1(\cA_{i-1, r+2})=q_*\bigl[\Bigl({g-1\choose
i-1}+c_1(\mathcal{F}) +
\frac{c_1^2(\mathcal{F})-2c_2(\mathcal{F})}{2}+\cdots \bigr)\cdot
\bigl(1-\frac{f^*c_1(\omega_p)}{2}+\frac{f^*c_1^2(\omega_p)}{12})+\ldots
\Bigr)\bigr].$$ Using that
$c_1(\mathcal{M})=p^*(\lambda)-c_1(\omega_p)$, one obtains that
$$c_1(\mathcal{F})={g-2\choose i-2} q^*(\lambda)+\Bigl((r+2){g-1\choose
i-1}-{g-2\choose i-2}\Bigr)f^*c_1(\omega_p)-{g-1 \choose i-1}
c_1(\mathfrak D).$$ We also use the identity
$$c_2(\mathcal{F})=c_2(f^*\wedge^{i-1}\mathcal{M})+\Bigl({g-1\choose
i-1}-1\bigr)c_1(f^*\wedge^{i-1}\mathcal{M})\cdot
\bigl((r+2)f^*c_1(\omega_p)-c_1(\mathfrak D)\bigr)+$$
$$+\frac{1}{2}{g-1\choose i-1}\Bigr({g-1\choose i-1}-1\Bigr)
\bigl((r+2)f^*c_1(\omega_p)-c_1(\mathfrak D)\bigr)^2,$$ which
together with the formula
$c_2(\mathcal{M})=c_1^2(\omega_p)-c_1(\omega_p)\cdot p^*(\lambda)$
and Lemma \ref{rules}, enable us to compute $c_1(\cA_{i-1, r+2})$.
In a similar fashion, we obtain from Grothendieck-Riemann-Roch
applied to the map $q$, that
$$c_1(\cA_{0, r+1})=(6r^2+6r+1)\lambda-(r+1)\sum_{j=1}^n \psi_j$$
and finally $$c_1(\pi^*\wedge^{i} \mathbb E\otimes \cA_{0, r+1})={g
\choose i}c_1(\cA_{0, r+1})+((2r+1)(g-1)-r){g-1 \choose i-1}\lambda
 ,$$
 which quickly leads to the stated formula.
\end{proof}

To compute the remaining
coefficients in $[\overline{\mathfrak{Mrc}}_{g, i}^r]$ we extend the
vector bundles  $\cA_{a, b}$ to sheaves over $\mm_{g, n}$ as
follows. We denote by $q:\mm_{g, n+1}\rightarrow \mm_{g, n}$ the
projection dropping the $(n+1)$-st marked point and by $\pi:\mm_{g, n}\rightarrow \mm_g$ the forgetful map. We introduce the
following twist of the Hodge bundle on $\mm_{g, n}$:
$$\H:=q_*\Bigl(\omega_q\otimes \OO_{\mm_{g, n+1}}\bigl(\sum_{[g/2]\leq j\leq g-1} \sum_{|S|\leq n}
(g-j)\Delta_{j:S}\bigr)\Bigr).$$ (In other words, the fibre of
$\H$ over a pointed curve from $\pi^*(\Delta_j)$ where $[g/2]\leq
j\leq g$, is the space of global sections of the genus $j$-aspect of
the limit $\mathfrak g^{g-1}_{2g-2}$ corresponding to the canonical
linear series). We then define
$$\cM:=\mbox{Ker}\{q^*(\H)\rightarrow \OO_{\mm_{g,
n+1}}\bigl(\sum_{[g/2]\leq j\leq g-1}\sum_{|S|\leq n}
(g-j)\Delta_{j: s}\bigr)\}.$$
Furthermore, for each pair of
integers $0\leq a\leq i, b\geq r+2$ or $(a, b)=(0, r+1)$, we define
$$\cA_{a, b}:=q_*\Bigl(\wedge^a \cM\otimes \omega_q^{\otimes
b}\otimes \OO_{\mm_{g, n+1}}\bigl(\sum_{[g/2]\leq j\leq
g-1}\sum_{|S|\leq n} ((2b-1)(g-j)-b)\Delta_{j:S} -\sum_{j=1}^n
\Delta_{0: j, n+1}\bigr)\Bigr).$$ (Obviously, this is an extension
of the definition of $\cA_{a, b}$ over $\cM_{g, n}$). The twists
were chosen in such a way that we have exact sequences of the type
$$0\longrightarrow \cA_{a, b}\longrightarrow \wedge^a \H\otimes
\cA_{0, b}\longrightarrow \cA_{a-1, b+1}\longrightarrow 0,$$ at
least in a dense open subset inside $\pi^{-1}(\cM_{g}\cup
\Delta_0\cup \Delta_1)$. Also, there exists a morphism $\phi:
\wedge^i \H\otimes \cA_{0, r+1}\rightarrow \cA_{i-1, r+2}$ which
over $\cM_{g, n}$ restricts to the vector bundle map defined in
Proposition \ref{hyper}.

\noindent \emph{Proof of Theorem \ref{mrc}.} We expand the class of
$\overline{\mathfrak{Mrc}}_{g, i}^r$ in $\mbox{Pic}(\mm_{g, n})$:
$$\overline{\mathfrak{Mrc}}_{g:i}^r\equiv
A\lambda+B\sum_{j=1}^n \psi_j -B_{irr}\delta_{irr}-\sum_{j, s\geq 0}
B_{j:s} \sum_{|S|=j} \delta_{j:S}.
$$ We have already determined the values of $A$ and $B$. One can write down the following relation
in $\mbox{Pic}(\mm_{g, n})$:
\begin{equation}\label{mrc6} c_1(\cA_{i-1, r+2}-\wedge^i \H\otimes \cA_{0,
r+1})=[\overline{\mathfrak{Mrc}}_{g, i}^r]+\sum_{j, s \geq 0} d_{j,
s}\sum_{|S|=s}\delta_{j:S}, \end{equation}
 where $d_{j, s}$ is the
multiplicity of the divisor $\Delta_{j:S}$ in the degeneracy locus
of $\phi$. By intersecting both sides of (\ref{mrc6}) with test
curves in $\mm_{g, n}$, sometimes  we are able to show that $\phi$
is generically non-degenerate along $\Delta_{j:S}$ (that is,
$d_{j:s}=0$), and then we explicitly determine the value of
$B_{j:S}$ in Theorem \ref{mrc}, otherwise we only get lower bounds
on $B_{j:S}$. We are only going to explain in detail the case of the
coefficient $B_{irr}$ the remaining ones being somewhat similar.

We define a test curve in the boundary of $\mm_{g, n}$ as follows.
If $[C, q, x_1, \ldots, x_n]\in \cM_{g-1, n+1}$ is a general pointed
curve, then we define the $1$-dimensional family $$C^0_n:=\{[C/y\sim
q, x_1, \ldots, x_n]\}_{y\in C}\subset \Delta_{irr}\subset \mm_{g,
n}.$$ The fibre of this family when the variable point $y\in C$ hits
the fixed marked point $x_i$ for $1\leq i\leq n$ is the pointed
curve $(\tilde{C}_{x_i}:=C\cup_{x_1, q} \PP^1, \tilde{x}_1, x_2,
\ldots, x_n)$, where $\tilde{x}_1\in \PP^1$ (here we regard $x_1,
\tilde{x}_1, q\in \PP^1$ as three distinct points). One has the
identities
$$C^0_n\cdot \delta_{irr}=-2g+2, \ C_n^0\cdot \delta_{1:\emptyset}=1, \ C_n^0\cdot \psi_i =1
\mbox{ for  } 1\leq i\leq n, \ C^0_n\cdot \lambda=C^0_n\cdot
\delta_{i:S}=0 \mbox{ for } (i, S)\neq (1, \emptyset).$$By
intersecting both sides of (\ref{mrc6}) with $C^0_n$ one can write
down the identity $C^0_n\cdot \overline{\mathfrak{Mrc}}_{g,
i}^r=(2g-2)B_{irr}+nB-B_{1:\emptyset}. $ On the other hand one also
has the relation $$A-12B_{irr}+B_{1:\emptyset}=0,$$ reflecting the
fact that $\overline{\mathfrak{Mrc}}_{g, i}^r$ is physically
disjoint from the curve $\{[C\cup_q R, x_1, x_2, \ldots, x_n]\}_{R}$
obtained by attaching to a fixed Brill-Noether general curve $[C, q,
x_1, \ldots, x_n]\in \cM_{g-1, n+1}$ a pencil of plane cubics in
which $R$ denotes a generic member. Thus determining $B_{irr}$ and
$B_{1:\emptyset}$ boils down to (i) showing that $\phi$ is
generically non-degenerate along $C^0_n$ and (ii) estimating the
intersection number $C^0_n\cdot c_1(\cA_{i-1, r+2}-\wedge^i
\H\otimes \cA_{0, r+1})$. By local analysis one can see that for
$1\leq l\leq i-1$ there are exact sequences of bundles over $C^0_n$
$$0\longrightarrow \cA_{i-l, r+l+1| C^0_n}\longrightarrow
\wedge^{i-l}\H\otimes \cA_{0, r+l+1 | C^0_n}\longrightarrow
\cA_{i-l-1, r+l+2 |C^0_n}\longrightarrow 0,$$ therefore we can write
the identities
$$C^0_n\cdot c_1(\cA_{i-1, r+2})=\sum_{l=1}^i
(-1)^{l-1}c_1\bigl(\wedge^{i-l}\H_{| C^0_n}\otimes \cA_{0, r+l+1 |
C^0_n}\bigr)=$$ $$ \sum_{l=1}^i (-1)^{l-1}\Bigl({g-1\choose
i-l-1}((2r+2l+1)(g-1)-n)c_1(\H_{| C^0_n})+{g \choose i-l}c_1(\cA_{0,
r+l+1 | C^0_n})\Bigr).$$ Next we describe the vector bundle $\cA_{0,
j | C^0_n}$. We identify $C^0_n$ with $C$ via the map $$C\ni
y\mapsto [C/y\sim q, x_1, \ldots, x_n]\in \mm_{g, n}$$ and denote by
$p_1, p_2:C\times C\rightarrow C$ the two projections, by
$\Delta\subset C\times C$ the diagonal and  set
$\Gamma_q:=\{q\}\times C\subset C\times C$. Then for every $j\geq r$
we have the following exact sequence of vector bundles on $C$:
$$0\longrightarrow (p_2)_*\Bigl(p_1^*(K_C^{\otimes j})\otimes
\OO_{C\times C}\bigl((j-1)\Delta+(j-1)\Gamma_q-\sum_{j=1}^n
\{x_j\}\times C\bigr)\Bigr)\longrightarrow \cA_{0, j |
C^0_n}\longrightarrow
$$
$$\longrightarrow (p_2)_*\Bigl(p_1^*(K_C^{\otimes j})\otimes
\OO_{C\times C}\bigl(j\Delta+j\Gamma_q-\sum_{j=1}^n \{x_j\}\times
C\bigr)\otimes \OO_{\Gamma_q}\Bigr)\otimes \OO_C(-x_1-\cdots
-x_n)\longrightarrow 0,$$ which quickly leads to the formula
$$c_1(\cA_{0, j| C^0_n})=1+2j-2jrg-2jg-j^2+j^2g+2jr+2ji.$$
Since one also has that $c_1(\H_{| C^0_n})=1$ (use that
$\H(y)=H^0(K_C\otimes \OO_C(y+q))$ for each $y\in C^0_n$ under the
identification described above), we obtain a formula for $C^0_n\cdot
\overline{\mathfrak{Mrc}}_{g, i}^r$ and ultimately a formula for
$B_{irr}$. Dealing with the other coefficients $B_{j:S}$ is similar
in general.

\hfill $\Box$

\subsection{Divisors defined by imposing linear conditions on marked points}\hfill

Here we present another general construction that produces families
of effective divisors on $\mm_{g, n}$. Like in Section 2, we pick
integers $g, r, d\geq 1$ such that $\rho(g, r, d)=0$, therefore we
can write $d=rs+r$ and $g=rs+s$ for some integer $s\geq 1$. We set
$n:=r+1$ and define the following divisor on $\cM_{g, n}$:
$$\mathfrak{Lin}^r_d:=\{[C, x_1, \ldots, x_{r+1}]: \exists L\in
W^r_d(C) \mbox{ such that } h^0(L\otimes \OO_C(-x_1-\cdots
-x_{r+1}))\geq 1\}.$$ We recall that we have introduced the number
$$N:=c_r=g!\frac{1! \ 2!\  \ldots r!}{(g-d+r)!\  \ldots (g-d+2r)!},$$
which counts  linear series $\mathfrak g^r_d$ on a general curve of
genus $g$ (cf. \cite{ACGH}). Our main result is the computation of
the class $[\overline{\mathfrak{Lin}}^r_d]$:
\begin{theorem}\label{pointedbn}
Fix integers $r, s\geq 1$ and set $d:=rs+r, g:=rs+s$. Then
$\overline{\mathfrak{Lin}}^r_d$ is an effective divisor on $\mm_{g,
n}$ and we have the following formula for its class in
$\rm{Pic}$$(\mm_{g, r+1})$:
$$\overline{\mathfrak{Lin}}^r_d\equiv \frac{r c_r}{rs+s-1}\bigl(a\lambda+c\sum_{j=1}^{r+1}
\psi_j-b_{irr}\delta_{irr}-\sum_{j, t\geq 0}
b_{j:t}\sum_{|S|=t}\delta_{j:S}\bigr),$$ where
$$a=\frac{(r+2)(r^2s^3-r^2s+2rs^3+6rs^2-2rs-8r+s^3+6s^2+3s-8)}{2(s+r+1)(rs+s-2)},$$
$$c=\frac{s+1}{2} , \mbox{ } b_{irr}=\frac{(s-1)(s+1)(r+1)(r+2)(rs+s+4)}{12(s+r+1)(rs+s-2)} ,$$
$$b_{j:0}=\frac{j(r+2)\bigl(rs(s^2-1)(r+2)+s(s^2-2j-3)+(r+1)(3s^2-js^2+2j-2)\bigr)}{2(r+s+1)(rs+s-2)}\ \mbox{ for }j\geq 1$$
$$b_{0:t}=\frac{t}{2r}(trs+ts-t+r-s+1) \mbox{ for } 2\leq
t\leq r+1,$$
$$b_{1:t}={t-1\choose 2}\frac{rs+s-1}{r}+\frac{(s-1)(s+1)(r+1)(r^3s+3r^2s-2s+4)}{2r(r+s+1)(rs+s-2)} \mbox{for all }
 t\geq 1$$
 $$\mbox{ and } b_{j: t}\geq b_{0: t} \mbox{ for all } j\geq 2.$$
\end{theorem}

\begin{remark} For $s=1$ and $r=g-1$, Theorem \ref{pointedbn} specializes to Logan's formula for the class of the
divisor $\overline{\mathfrak{Lin}}^{g-1}_{2g-2}$ of points $[C, x_1,
\ldots, x_g]\in \cM_{g, g}$ with $h^0(\OO_C(x_1+\cdots +x_g))\geq
2$. We have the formula (cf. \cite{Log}, Theorem 5.4):
$$\overline{\mathfrak{Lin}}^{g-1}_{2g-2}\equiv -\lambda+\sum_{j=1}^g
\psi_j-0\cdot \delta_{irr}-\sum_{t= 2}^{g} {t+1\choose
2}\sum_{|S|=t}\delta_{0:S}- \sum_{t=1}^{g} {t\choose 2} \sum_{|S|=t} \delta_{1:S}- \cdots.$$
In the next case, $s=2, g=2r+2$ and $d=3r$ we get a new divisor on $\mm_{2r+2, r+1}$ and our formula reads
\begin{equation}\label{pointed2}
\overline{\mathfrak{Lin}}^r_{3r}\equiv
\frac{1}{2(2r+1)}\Bigl((3r+5)(r+2)\lambda+3r\sum_{j=1}^{r+1}
\psi_j-{r+2\choose 2}\delta_{irr}-  \sum_{|S|=2}
\delta_{0:S}-\cdots\Bigr).
\end{equation}
\end{remark}

\noindent \emph{Proof of Theorem \ref{pointedbn}.} Proving that
$\mathfrak{Lin}^r_d$ is a divisor on $\cM_{g, r+1}$ follows
immediately from Brill-Noether theory: a general $[C]\in \cM_g$ has
precisely $N$ linear series $L\in W^r_d(C)$ and each of them is base
point free and satisfies $h^0(L)=r+1$ and $\mathfrak{Lin}^r_d$
consists of those $(r+1)$-tuples of points on $C$ which are not in
general linear position with respect to some $L\in W^r_d(C)$.  We
start the calculation of the class of
$\overline{\mathfrak{Lin}}^r_d$ by determining the coefficients of
$\lambda, \delta_{irr}$ and $\psi_j$. To do this we use the
observation that if $\pi_n:\mm_{g, n}\rightarrow \mm_{g, n-1}$ is
the projection forgetting the marked point labeled by $n$ and $D$ is
any divisor class on $\mm_{g, n}$, then for distinct labels $i, j
\neq n$, the $\lambda$, $\delta_{irr}$ and $\psi_j$ coefficients of
the divisors $D$ on $\mm_{g, n}$ and $(\pi_n)_*(D\cdot \delta_{0:i
n})$ on $\mm_{g, n-1}$ are the same (see \cite{FMP}, Prop. 4.4). The
divisor $(\pi_n)_*(D\cdot \delta_{0:i n})$ can be thought of as the
locus of those points $[C, x_1, \ldots, x_n]\in D$ where the points
$x_i$ and $x_n$ are allowed to come together. Using this observation
repeatedly, we obtain that the divisor
$\overline{\mathfrak{Lin}}_d^r(1)$ on $\mm_{g, 1}$ obtained by
letting all points $x_1, \ldots, x_n$ come together, has the same
$\lambda$ and $\delta_{irr}$ coefficients as $\mathfrak{Lin}^r_d$.
But clearly
$$\mathfrak{Lin}_d^r(1)=\{[C, x]\in \cM_{g, 1}: \exists L\in W^r_d(C) \mbox{ such
that } h^0(L\otimes \OO_C(-(r+1)x))\geq 1\},$$ that is,
$\mathfrak{Lin}_d^r(1)$ is generically the locus of ramification
points in one of the finitely many linear series $\mathfrak g^r_d$
on a given curve of genus $g$. By applying Theorem 4.1 from
\cite{EH2}, we obtain that the class of
$\overline{\mathfrak{Lin}}_d^r(1)$ can be written as a combination
$\overline{\mathfrak{Lin}}_d^r(1)\equiv \mu \cdot BN+\nu \cdot
\mathcal{W}$, where
$$BN:=(g+3)\lambda-\frac{g+1}{6} \delta_{irr}-\sum_{j=1}^{g-1}
j(g-j)\delta_{j: 1}$$ is the pull-back from $\mm_g$ of the class of
the Brill-Noether divisor and
$$\mathcal{W}:=-\lambda+{g+1\choose 2}\psi-\sum_{j=1}^{g-1}
{g-j+1\choose 2} \delta_{j:1}$$ is the class of the Weierstrass
divisor. To determine the coefficients $\mu$ and $\nu$ we use two
test curves inside $\mm_{g, 1}$. First we fix a genus $g$ curve $C$
and we let the marked point vary along $C$. If we denote this curve
by $\bar{C}\subset \mm_{g, 1}$, then the only generator of
$\mbox{Pic}(\mm_{g, 1})$ which has non-zero intersection number with
$\bar{C}$ is $\psi$, and $\bar{C}\cdot \psi=2g-2$. On the other hand
$\bar{C}\cdot \overline{\mathfrak{Lin}}_d^r(1)$ is the total number
of ramification points on all $\mathfrak g^r_d$'s on $C$. This
number is $N(r+1)(d+r(g-1))$ (see e.g. \cite{EH1}, pg. 345), which
shows that $\nu=N(r+1)(d+r(g-1))/((g-1)g(g+1))$. To compute $\mu$ we
use a second test curve constructed as follows: we fix a $2$-pointed
elliptic curve $[E, x, y]\in \cM_{1, 2}$ such that the class $x-y\in
\mbox{Pic}^0(E)$ is not torsion, and a fixed general curve $[C]\in
\cM_{g-1}$. We define the family $\bar{C}_1:=\{(C\cup _y E,
x)\}_{y\in C}$ (that is, the point of attachment varies on $C$). The
only generator of $\mbox{Pic}(\mm_{g, 1})$ meeting $\bar{C}_1$
non-trivially is $\delta_{1:1}=\delta_{g-1: \emptyset}$, in which
case $\bar{C}_1\cdot \delta_{1;1}=-2g+4$. The calculation of
$\bar{C}_1\cdot \overline{\mathfrak{Lin}}_d^r(1)$ is a standard
exercise in the theory of limit linear series. Suppose $l=\{l_C,
l_E\}$ is a limit $\mathfrak g^r_d$ on $C\cup _y E$ such that
$a_r^{l_E}(x)\geq r+1$. Then because the class $x-y \in
\mbox{Pic}^0(E)$ is not torsion, we must have that $\rho(l_E, x,
y)=0$ and $\rho(l_C, y)=-1$. An easy calculation shows that we must
also have $a^{l_E}(x)=(0, 1, \ldots, r-1, r+1)$ and $a^{l_C}(y)=(0,
2, 3, \ldots, r, r+2)$, and moreover the aspect $l_E$ is uniquely
determined. Thus we have to count the number of points $y\in C$ such
that there exists $L\in W^r_d(C)$ with the property that
$h^0(L(-2y))\geq r$ and $h^0(L(-(r+2)y))\geq 1$.

To compute this number we further degenerate the curve $C$ to a
transverse union $R\cup _{y_1} E_1\cup \ldots \cup_{y_{g-1}}
E_{g-1}$ consisting of a smooth rational spine $R$ and $g-1$
elliptic tails $E_1, \ldots, E_{g-1}$. Using Proposition 1.1 from
\cite{EH1} we see that the point $y$ has to specialize to one of the
tails $E_j$, and without loss of generality we assume that $y\in
E_1$ (all the intersection numbers we compute will be multiplied by
$g-1$ to account for $y$ lying on a different elliptic tail).
Suppose now that $l=\{l_R, l_{E_1}, \ldots, l_{E_{g-1}}\}$ is a
limit $\mathfrak g^r_d$ on $R\cup E_1\cup \ldots \cup E_{g-1}$ such
that $a^{l_{E_1}}(y)=(0, 2, 3, \ldots, r, r+2)$. Then $\rho(l_R,
y_1, \ldots, y_{g-1})=0, \rho(l_{E_j}, x_j)=0$ for $2\leq j\leq g-1$
and $\rho(l_{E_1}, y_1, y)=-1$. A close inspection shows that there
are three numerical possibilities: \vskip 2pt \noindent
{\textbf{($\alpha$)}} $a^{l_R}(y_1)=(0, 2, 4, 5, \ldots, r, r+1,
r+3)$ and then $y_1-y\in \mbox{Pic}^0(E_1)[2]$. This contribution
will be equal to $3(g-1)$ multiplied by the number of $\mathfrak
g^r_d$'s on $R$ having ordinary cusps at $g-2$ general points and
vanishing $(0, 2, 4, \ldots, r+1, r+3)$ at another fixed point. By
Schubert calculus this number equals the product of Schubert cycles
$\sigma_{(0, 1, \ldots, 1)}^{g-2}\cdot \sigma_{(0, 1, 2, \ldots, 2,
3)} \in H^*(\GG(r, d))$. \vskip 1pt
 \noindent {\textbf{($\beta $)}} $a^{l_R}(y_1)=(0, 3, 4, \ldots,
r+1, r+2)$, in which case $y-y_1\in  \mbox{Pic}^0(E_1)[r+2]$. The
number we get in this situation is $(g-1)((r+2)^2-1)\sigma_{(0, 1,
\ldots, 1)}^{g-2} \cdot \sigma_{(0, 2, \ldots , 2)}\in H^*(\GG(r,
d))$. \vskip 1pt \noindent {\textbf{($\gamma$)}} $a^{l_r}(y)=(1, 2,
4, \ldots, r+1, r+2)$ and then $y-y_1\in \mbox{Pic}^0(E_1)[r]$. We
obtain a final contribution of $(g-1)(r^2-1)\sigma_{(0, 1, \ldots,
1)}^{g-2}\cdot \sigma_{(0, 0, 1, \ldots, 1)}\in H^*(\GG(r, d-1))$.
Adding all these together and using (\ref{schubert0}), we obtain
that the total intersection number is
$$\bar{C}_1\cdot
\overline{\mathfrak{Lin}}_d^r(1)=\frac{N\
r(r+1)(r+2)(rs+2s^2-4+s)}{s+r+1},$$ which leads to
$$\mu=\frac{N r(r+1)(r+2)(s-1)(s+1)(rs+s+4)}{2(s+r+1)(rs+s-2)(rs+s-1)(rs+s+1)}$$
and then the stated formulas for the $\lambda$ and $\delta_{irr}$
coefficients. We also note that the $\delta_{j: \emptyset}$
coefficient of $\overline{\mathfrak{Lin}}^r_d$ equals the
$\delta_{j:\emptyset}=\delta_{g-j:1}$ coefficient of
$\overline{\mathfrak{Lin}}^r_d(1)$ and this is equal to
$j(g-j)\mu+j(j+1)\nu/2$ and we obtain the desired expression for
$b_{j:0}$. Next we determine the coefficient $c$ of the $\psi_j$
classes. We introduce the divisor $\overline{\mathfrak{Lin}}_d^r(2)$
on $\mm_{g, 2}$ obtained by letting $x_2, \ldots, x_{r+1}$ coincide
while keeping $x_1$ apart. Then $$\mathfrak{Lin}^r_d(2)=\{[C, x_1,
x_2]\in \cM_{g, 2}: \exists L\in W^r_d(C) \mbox{ such that }
h^0(L\otimes \OO_C(-x_1-rx_2))\geq 1\}$$ and we can write
$\overline{\mathfrak{Lin}}^r_d(2)\equiv N(a \lambda+c \psi_1+c_2
\psi_2 -e_{12} \delta_{0:12}-\cdots)$ (that is, the $\lambda$ and
$\psi_1$ coefficients coincide with those of
$\overline{\mathfrak{Lin}}^r_d$). We intersect
$\overline{\mathfrak{Lin}}^r_d(2)$ with two curves in $\mm_{g, 2}$:
consider a curve $C$ of genus $g$ and define $\tilde{C}_1:=\{[C,
x_1, x_2]\}_{\{x_2 \ moves \ on \ C \}}$ and $\tilde{C}_2:=\{[C,
x_1, x_2]\}_{\{x_1 \  moves \ on \ C \}}$. Then $\tilde{C}_2\cdot
\overline{\mathfrak{Lin}}^r_d(2)=N(c_2+(2g-1)c-e_{12})=N(d-r)$ and
$\tilde{C}_1\cdot
\overline{\mathfrak{Lin}}^r_d(2)=N((2g-1)c_2+c-e_{12})=N(r+1)(d-1+(r-1)(g-1))$.
(The first identity is obvious, for the second, use that what we are
counting is the total number of ramification points on all linear
series $L\otimes \OO_C(-x_1)$, where $L\in W^r_d(C)$ and $x_1\in C$
is a fixed general point). We thus have a system of two equations in
the unknowns $c, c_2$ and $e_{12}$, but we can also use that
$e_{12}$ equals the $\psi$ coefficient of
$\overline{\mathfrak{Lin}}^r_d(1)=(\pi_2)_*(\overline{\mathfrak{Lin}}^r_d(2)\cdot
\delta_{0: 12})$, where $\pi_2:\mm_{g, 2}\rightarrow \mm_{g, 1}$ is
the map forgetting the second point. Thus $e_{12}=\nu g(g+1)/2$,
which gives us enough relations to determine $c$. We note that in
this way we also determine $b_{0:2}=(2rs+r+s-1)/(rs+s-1)$.

 To compute the
coefficient $b_{j: t}$ for $1\leq t\leq r+1$ we consider another
test curve defined as follows: we fix integers $1\leq j\leq g-1$ and
$1\leq t\leq r+1$, together with general pointed curves $[C, y, x_2,
\ldots, x_t]\in \cM_{j, t}$ and $[Y, y, x_{t+1}, \ldots, x_{r+1}]\in
\cM_{g-j, r-t+2}$. We define the test curve $\bar{C}_{j,
s}:=\{C\cup_y Y, x_1, \ldots, x_t, x_{t+1}, x_{r+1}\}_{x_1\in C}$
(thus $x_1$ is the moving point on the genus $j$ component). Then we
have the relation
\begin{equation}\label{recursion} \bar{C}_{j,t}\cdot
\overline{\mathfrak{Lin}}^r_d=(2j+2t-3)c-(t-1)b_{0:2}+b_{j:t}-b_{j:t-1},
\end{equation}
which can be used to compute $b_{j:t}$ provided we know $b_{j:t-1}$
(note that we have already computed $b_{j:\emptyset}$ for all $j$).

We now describe directly the intersection cycle $\bar{C}_{j, t}\cdot
\overline{\mathfrak{Lin}}^r_d$. Since $[C, y]\in \cM_{j, 1}$ and
$[Y, y]\in \cM_{g-j, 1}$ are general, on the stable $Y\cup_y C$
there will be precisely $N$ limit $\mathfrak g^r_d$'s which can be
described as follows: We first choose a Schubert (ramification)
sequence $\mbox{max}\{0, j-t\}\leq \alpha_0 \leq \ldots \alpha_r\leq
j$ such that $\sum_{i=0}^r \alpha_i=rj$. Then we choose $l_Y\in
G^r_d(Y)$ having vanishing sequence $(a^{l_Y}(y)=\alpha_i+i)_{0\leq
i\leq r}$; in fact there will be $\sigma_{(0, 1, \ldots,
1)}^{g-j}\cdot \sigma_{(\alpha_0, \ldots, \alpha_r)}\in H^*(\GG(r,
d))$ such linear series. On $C$ we choose a complementary linear
series $l_C\in G^r_d(C)$ with vanishing sequence
$(a_i^{l_C}(y)=rs+i-\alpha_i)_{0\leq i\leq r}$; there are
$\sigma_{(0, 1, \ldots, 1)}^j\cdot \sigma_{(rs-\alpha_r, \ldots,
 rs-\alpha_0)}\in H^*(\GG(r, d))$ choices. Every limit linear series
on $C\cup _y Y$ appears in this way and the intersection
$\bar{C}_{j, t}\cdot \overline{\mathfrak{Lin}}^r_d$ is everywhere
transverse (cf. \cite{EH3}). We also have the identity
\begin{equation}\label{schubert}
N=\sum_{j-t\leq \alpha_0\leq \ldots \leq \alpha_r\leq j,\ \sum
_{l=0}^r \alpha_l=j} \bigl(\sigma_{(0, 1, \ldots, 1)}^{g-j}\cdot
\sigma_{(\alpha_0, \ldots, \alpha_r)}\bigr) \ \bigl(\sigma_{(0, 1,
\ldots, 1)}^j\cdot \sigma_{(rs-\alpha_r, \ldots,
rs-\alpha_0)}\bigr).
\end{equation}
If $l=\{l_C, l_Y\}$ is one of these $N$ limit $\mathfrak g^r_d$'s
corresponding to a sequence $(\alpha_0\leq \ldots \leq \alpha_r)$ as
above, the condition that there exists $x_1\in C$ such that the
divisor $x_1+\cdots+x_{r+1}$ is the specialization of a linear
divisor with respect to a $\mathfrak g^r_d$ on a nearby smooth
curve, can be translated as follows: there exist sections $\sigma_Y
\in |l_Y|, \sigma_C \in |l_C|$ such that $\mbox{div}(\sigma_Y)\geq
x_{t+1}+\cdots+x_{r+1}$ and $\mbox{div}(\sigma_C)\geq
x_1+\ldots+x_t$; the sections $\sigma_C$ and $\sigma_Y$ being the
limit linear series specializations of a single section on a nearby
smooth curve, they must also satisfy the compatibility relation
$\mbox{ord}_y(\sigma_Y)+\mbox{ord}_y(\sigma_C)=rs+r$. Because the
fixed points $x_{t+1}, \ldots, x_{r+1}\in Y$ are general, they
impose independent conditions on $l_Y$ which quickly leads to the
equalities
$\mbox{ord}_y(\sigma_Y)=a^{l_Y}_{t-1}(y)=\alpha_{t-1}+t-1$, hence
$\mbox{ord}_y(\sigma_C)=a_{r-t+1}^{l_C}(y)$. Thus
$\mbox{div}(\sigma_C)\geq a_{r-t+1}^{l_C}(y)+x_2+\ldots+x_b$ and up
to multiplication by scalars, the sections $\sigma_C$ and $\sigma_Y$
are unique with this property. For each $\sigma_C$ we have precisely
$d-a_{r-t+1}^{l_C}(y)-(t-1)=\alpha_{t-1}$ choices for $x_1\in C$.
Therefore
\begin{equation}\label{schubert2}
\bar{C}_{j, t}\cdot \overline{\mathfrak{Lin}}^r_d=\sum_{\alpha_0\leq
\cdots \leq \alpha_r} \alpha_{t-1}\bigl(\sigma_{(0, 1, \ldots,
1)}^{g-j}\cdot \sigma_{(\alpha_0, \ldots, \alpha_r)}\bigr)\cdot
\bigl(\sigma_{(0, 1, \ldots, 1)}^j\cdot \sigma_{(rs-\alpha_r,
\ldots, rs-\alpha_0)}\bigr).
\end{equation}
For $j=0$ the only sequence $(\alpha_l)_{0\leq l \leq r}$ allowed is
the sequence $(0, \ldots, 0)$ which shows that $\bar{C}_{0:t}\cdot
\overline{\mathfrak{Lin}}^r_d=0$ for all $3\leq t\leq r+1$. Since
$b_{0:2}$ has already been determined, applying (\ref{recursion}) we
obtain the stated formulas for $b_{0:t}$. Similarly, for $j=1$ the
only sequence allowed is $(0, 1, \ldots, 1)$ and then
$\bar{C}_{1:t}\cdot \overline{\mathfrak{Lin}}^r_d=N$ for $t\geq 2$,
while $\bar{C}_{1:1}\cdot \overline{\mathfrak{Lin}}^r_d=0$; this
allows us to determine $b_{1:t}$ for all $t$. When $j\geq 2$ for
each sequence $(\alpha_l)_{0\leq l\leq r}$ appearing in this sum, we
have the inequalities $rj=\sum_{l=0}^r \alpha_l \leq
t\alpha_{t-1}+(r+1-t)j$, therefore $\alpha_{t-1}\geq (t-1)j/t$ and
then $\bar{C}_{j, t}\cdot \overline{\mathfrak{Lin}}^r_d\geq
N(t-1)j/t$. To obtain the desired bound on $b_{j:t}$ we use
repeatedly (\ref{recursion}) and we can write
$$b_{j:t}-b_{j:0}={t \choose 2}
b_{0:2}-(2bj+b^2-2b)c+\sum_{l=1}^t \bar{C}_{j, l}\cdot
\overline{\mathfrak{Lin}}^r_d.$$ Using the previous inequality we
can now check
 that $b_{j:t}\geq b_{j: 0}$.

\hfill $\Box$

\subsection{The divisor of $n$-fold points}
\hfill

We describe another way of constructing effective divisors on
$\mm_{g, n}$. Instead of looking at loci of points $[C, x_1, \ldots,
x_n]\in \cM_{g, n}$ for which the points $x_1, \ldots, x_n$ become
linearly dependent in a suitable embedding of $C$, we can consider
the loci where the marked points give rise to an $n$-fold point on a
suitable model of $C$. Given $[C]\in \cM_g$ and a linear series $l=(L, V)\in
G^r_d(C)$, we say that the divisor $\Gamma:=x_1+\cdots+x_n$ is an
\emph{$n$-fold point} for $C$ and $l$ if $\mbox{dim}\bigl(V\cap
H^0(L\otimes \OO_C(-\Gamma))\bigr)\geq r.$
\begin{definition}
Fix integers $g, r, d, n\geq 1$ such that $\rho(g, r, d)-r(n-1)=-1$.
We define the locus of $n$-fold points in $\cM_{g, n}$
$$\mathfrak {Nfold}_{g, d}^r:=\{[C, x_1, \ldots, x_n]\in \cM_{g,
n}: \exists L\in W^r_d(C) \mbox{ with } \mbox{dim }
H^0(L(-x_1-\cdots -x_n))=r\}.$$
\end{definition}

We have computed the class of $\mathfrak{Nfold}_{g, d}^r$ in the
case $r=1$. The calculation is along the same lines as that of the
class of $\overline{\mathfrak{Lin}}^r_d$ in Theorem \ref{pointedbn}:
 \begin{theorem}\label{misc}
  Fix integers $g\geq 1$ and
$n\geq 0$ such that $d:=(g+n)/2\in \mathbb Z$. The class of the
compactification of the divisor $\mathfrak{Nfold}_{g, d}^1$ of $n$-fold points on $\cM_{g,
n}$ is given by the formula:
$$\overline{\mathfrak{Nfold}}_{g, d}^1\equiv \Bigl(\frac{10n}{g-2}{g-2
\choose d-1}-\frac{n}{g}{g\choose
d}\Bigr)\lambda+\frac{n-1}{g-1}{g-1\choose d-1}\sum_{j=1}^n
\psi_j-\frac{n}{g-2}{g-2\choose d-1}\delta_{irr}-$$ $$-\sum_{t\geq
2}\frac{t(n^2-g+tgn-tn)}{2(g-1)(g-d)}{g-1 \choose
d}\sum_{|S|=t}\delta_{0:S}-\cdots.
$$
\end{theorem}
\begin{proof} The coefficients of $\lambda, \delta_{irr}$ and
$\psi_j$ $(1\leq j\leq n)$ in the expansion of
$\overline{\mathfrak{Nfold}}_{g, d}^1$ equal the coefficients of $\lambda,
\delta_{irr}$ and $\psi$  respectively in the expansion
of the divisor $\overline{\mathfrak{Nfold}}_{g, d}^1(1)$ on $\mm_{g,
1}$ obtained from $\overline{\mathfrak{Nfold}}_{g, d}^1$ by letting
the points $x_1, \ldots, x_n\in C$ coalesce. Clearly,
$$\mathfrak{Nfold}_{g, d}^1(1):=\{[C, x]\in \cM_{g, 1}: \exists L\in
W^1_d(C) \mbox{ with } h^0(C, L(-n\cdot x))\geq 1 \},$$ and this is
a "pointed" Brill-Noether divisor on $\cM_{g, 1}$ in the sense of
\cite{EH2}. To compute the class of its compactification in $\mm_{g, 1}$ once again we use \cite{EH2},
Theorem 4.1 and  write $\overline{\mathfrak{Nfold}}_{g,
d}^1(1)\equiv \mu\cdot BN+\nu\cdot \mathcal{W}$, where the divisor
classes $BN$ and $\mathcal{W}$ have the same significance as in the
proof of Theorem \ref{pointedbn}. By applying \cite{Log}, Theorem
4.5, we find that
$$\mu=\frac{6n}{(g+1)(g-2)}{g-2 \choose d-1}\ \mbox{ and } \
\nu=\frac{n(n-1)(n+1)}{g(g-1)(g+1)}{g\choose d} .$$ The remaining
coefficients of $[\overline{\mathfrak{Nfold}}_{g, d}^1]$ are
determined by intersecting the locus
$\overline{\mathfrak{Nfold}}_{g, d}^1$ with the fibral test curves
lying entirely in the boundary divisors of $\mm_{g, n}$. The calculation is
straightforward and relies on Section 3 from \cite{Log}. We skip
these details.
\end{proof}

\section{The Kodaira dimension of $\mm_{g, n}$}

In this section we prove Theorem \ref{mgn}. We treat each case
individually but for each $g$ we only work out the case of the
minimal $n=n(g)$ for which our methods show that $\mm_{g, n(g)}$ is
of general type. From this it follows automatically  that $\mm_{g,
n}$ is of general type for all $g\geq g(n)$ (see \cite{Log}, Theorem
2.4).

\noindent {\emph{Proof of Theorem \ref{mgn}.}}
\noindent\textbf{[$\mm_{4, 16}$] and [$\mm_{6, 16}$]}. \mbox{ } We
consider the divisor $\mathfrak{Mrc}_{4, 0}^2$ on $\mm_{4, 15}$ introduced in
Theorem \ref{mrc}. We have seen that $\overline{\mathfrak{Mrc}}_{4, 0}^2\equiv
-37\lambda+3\sum_{j=1}^15\psi_j+3\delta_{irr}-7\sum_{|S|=2}\delta_{0:S}-\cdots.$
We consider the  maps $\pi_j:\mm_{4, 16}\rightarrow \mm_{4, 15}$
obtained by forgetting the marked point labeled by $1\leq j\leq 16$.
Then there exists a constant $\alpha>0$ such that
$$\sum_{j=1}^{16} (\pi_j)^*(\overline{\mathfrak{Mrc}}_{4, 0}^2)\equiv \alpha(-37\lambda+\frac{45}{16} \sum_{j=1}^16 \psi_j+
3\delta_{irr}-\frac{13}{2}\sum_{|S|=2}\delta_{0:S}-\cdots).$$ The
class of the Petri divisor on $\mm_4$ being (up to a $>0$ constant)
$17\lambda-2\delta_{irr}-\cdots$, we obtain that $K_{\mm_{4, 16}}$ is big,
being a positive combination of $\sum_{j=1}^{16}
(\pi_j)^*(\overline{\mathfrak{Mrc}}_{4: 0}^2)$, the pull-back of the Petri class,
an ample class and boundary divisors. The same argument works in the
case of $\mm_{6, 16}$ except that we start with the divisor
$\overline{\mathfrak{Mrc}}_{6, 0}^1$ on $\mm_{6, 15}$ which is pulled back
to $\mm_{6, 16}$ in all possible ways.

\noindent {\textbf{[$\mm_{5, 15}$]}}. \mbox{ } On $\mm_{5, 12}$ we
have the identity of divisor classes $$\overline{\mathfrak{Mrc}}_{5, 0}^1 \equiv
-13\lambda+2\sum_{j=1}^12\psi_j+\delta_{irr}-5\sum_{|S|=2}\delta_{0:S}-\cdots.$$
Pulling this class back to $\mm_{5, 15}$
 in all possible ways by forgetting sets of three marked points, we obtain that
 the class $-66\lambda+\frac{42}{5}\sum_{j=1}^{15} \psi_j+5\delta_{irr}-\cdots$ is
 effective on $\mm_{5, 15}$. Using the Brill-Noether class $8\lambda-\delta_{irr}- \cdots$ on
 $\mm_5$,
 we get that $K_{\mm_{5, 15}}$ is a big class.

\noindent {\textbf{[$\mm_{18, 9}$]}}. \ We use our divisor
$\overline{\mathfrak{Lin}}^8_{24}$: there is a positive constant
$\alpha$ such that $$\alpha \overline{\mathfrak{Lin}}^8_{24}\equiv
290\lambda +24 \sum_{j=1}^9 \psi_j-45\delta_{irr}-\cdots.$$
 On the
other hand the class of the multiple of the Petri divisor
$\overline{\mathcal{GP}}^{8}_{18, 24}$ on $\mm_{18}$ is equal to
$\frac{302}{45}\lambda-\delta_0-\sum_{j=1}^9 b_j \delta_j$, where
$b_j> 1$ for $j\geq 1$. It follows that we can write $K_{\mm_{18,
9}}$ as a positive combination of multiples of
 $\overline{\mathfrak{Lin}}^8_{24}$, $\pi^*(\overline{\mathcal{GP}}^8_{18, 24})$,
 boundary divisors and an ample class on $\mm_{18, 9}$.

\noindent {\textbf{[$\mm_{19, 7}$] and [$\mm_{14, 10}]$}}. \ In
these cases we use the divisors of $n$-fold points,
$\overline{\mathfrak{Nfold}}^1_{19, 13}$ on $\mm_{19, 7}$ and
$\overline{\mathfrak{Nfold}}^1_{14, 12}$ on $\mm_{14, 10}$
respectively. Using Theorem \ref{misc} we see
 that the canonical bundle of $\mm_{g, n}$ can be written as a positive combination of these divisors,
 the pull-back of the Brill-Noether divisor from $\mm_{19}$ and $\mm_{14}$ respectively, a suitable
  ample class and boundary divisors.

\noindent {\textbf{[$\mm_{15, 10}$]}. \mbox{ } We use a slightly
different technique. On $\mm_{15, 11}$ we have the divisor
$\overline{\mathfrak{Nfold}}^1_{15, 13}$ of points $[C, x_1, \ldots,
x_{11}]$ such that $x_1+\cdots +x_{11}$ appears in a fibre of a
$\mathfrak g^1_{13}$ on $C$. We push this divisor down to $\mm_{15,
10}$ by letting two of the points $x_j\in C$ coalesce, that is, we
define $$E:=\frac{1}{11}\sum_{j=1}^{10}
(\pi_j)_*(\overline{\mathfrak{Nfold}}_{15, 13}^1\cdot \delta_{0:j,
11}),$$ where $\pi_j:\mm_{15, 11}\rightarrow \mm_{15, 10}$ forgets
the marked point labeled by $j$. It is easy to check using Theorem
\ref{misc} that
$$E\equiv 33\lambda + \frac{396}{5} \sum_{j=1}^{10} \psi_j
-11\delta_{irr}-\cdots \bigl(\mbox{ use that }(\pi_j)_*(\psi_j\cdot
\delta_{0:j, 11})=(\pi_j)_*(\psi_{11}\cdot \delta_{0: j, 11})=0
\mbox{ } \bigr).$$ It turns out that $K_{\mm_{15, 11}}$ is in the
span of $E, \pi^*(\mm_{15, 14}^3)$, an ample class and boundaries.

\noindent {\textbf{[$\mm_{20, 6}$]}. From  \cite{Log} Theorem 5.4,
one knows that the class \ $-\lambda+\frac{22}{3}\sum_{j=1}^6
\psi_j-0\cdot \delta_{irr}-\cdots$ is effective on $\mm_{20, 6}$.
Next, if $\chi_{i, j}:\mm_{20, 6}\rightarrow \mm_{21}$ denotes the
map which associates to a $6$-pointed curve of genus $20$ a nodal
curve of genus $21$ obtained by identifying the marked points
labeled $i$ and $j$, we also get that the class $$\sum_{i<j}
\chi_{i, j}^*(\overline{\mathcal{Z}}_{21, 0})\equiv
c(\frac{2459}{377} \lambda+\frac{1}{3}\sum_{j=1}^6
\psi_j-\delta_{irr}-\cdots),$$ with $c>0$, is also effective on
$\mm_{20, 6}$. The conclusion now follows easily.

\end{document}